\documentclass[11pt]{amsart}
\usepackage{Macros}

\title[Cyclic sieving, necklaces, branching rules, and Thrall's problem]{Cyclic sieving, necklaces, and branching rules related to Thrall's problem}
\author{Connor Ahlbach and Joshua P.\ Swanson}
\address{Department of Mathematics, University of Washington,
Seattle, WA 98195, USA; Connor Ahlbach: ahlbach@uw.edu; Joshua P. Swanson: jps314@uw.edu}
\thanks{The authors were partially supported by the National Science
Foundation grant DMS-1101017.}
\date{\today}

\begin{document}

\begin{abstract}
  We show that the cyclic sieving phenomenon of
  Reiner--Stanton--White together with necklace generating
  functions arising from work of Klyachko offer a
  remarkably unified, direct, and largely bijective approach
  to a series of results due to Kra{\'s}kiewicz--Weyman,
  Stembridge, and Schocker related to the so-called
  higher Lie modules and branching rules for inclusions
  $ C_a \wr S_b \hookrightarrow S_{ab} $.
  Extending the approach gives monomial expansions for certain graded
  Frobenius series arising from a generalization of Thrall's
  problem.
\end{abstract}
\maketitle

\section{Introduction}\label{sec:intro}
The \textit{Lie module} $\cL_n$ is the $n$th degree component of the
free Lie algebra over $\bC$ with $m$ generators, which is naturally a
$\GL(\bC^m)$-module. The Lie modules were famously studied by
Thrall \cite{MR0006149} in the 1940's and have been extensively
studied by Brandt \cite{MR0011305}, Klyachko \cite{klyachko74},
Kra{\'s}kiewicz--Weyman \cite{MR1867283},
Garsia \cite{MR1039352}, Gessel--Reutenauer \cite{MR1245159},
Reutenauer \cite{MR1231799}, Sundaram \cite{MR1273390},
Schocker \cite{MR1984625}, and many others. Thrall more generally
introduced a certain $\GL(\bC^m)$-decomposition
$\oplus_{\lambda \in \Par} \cL_\lambda$ of the tensor algebra
of $\bC^m$ arising from the Poincar\'e--Birkhoff--Witt
theorem, where $\cL_{(n)} = \cL_n$. The
$\cL_\lambda$ are sometimes called the \textit{higher Lie modules}.
Thrall's original paper considered the determination of the
multiplicity of the irreducible $V^\mu$ in $\cL_\lambda$,
which is often referred to as \textit{Thrall's problem}.
This problem is still open 75 years later. See
\Cref{ssec:background_thrall} and \cite{MR1231799} for more
background on Thrall's problem and \cite{reiner15} for a recent
summary of related work. See \Cref{sec:background} for missing
definitions.

Kra{\'s}kiewicz--Weyman \cite{MR1867283} gave a combinatorial
solution to Thrall's problem when $\lambda = (n)$. In particular, they
showed the multiplicity of $V^\mu$ in $\cL_{(n)}$ is
  \[ \#\{T \in \SYT(\mu) : \maj(T) \equiv_n 1\}, \]
i.e.~the number of
standard tableaux of shape $\mu$ with major index $1$ modulo $n$. Their
argument crucially hinges upon the formula
\begin{equation}\label{eq:SYT_evals}
  \SYT(\mu)^{\maj}(\omega_n^r) = \chi^\mu(\sigma_n^r)
\end{equation}\label{eq:kw_csp}
where we write the major index generating function as
  \[ \SYT(\mu)^{\maj}(q) \coloneqq \sum_{T \in \SYT(\mu)}
     q^{\maj(T)}, \]
$\omega_n$ is a primitive $n$th complex root of unity, $\sigma_n$ is
an $n$-cycle in the symmetric group $S_n$, and $\chi^\mu$ is
the character of the $S_n$-irreducible indexed by a partition
$\mu$ of $n$. The analysis in \cite{MR1867283} is somewhat indirect.
It involves results of Lusztig and Stanley on coinvariant algebras
and an intricate though beautiful argument involving $\ell$-decomposable
partitions.

Equation~\eqref{eq:kw_csp} bears a striking resemblance
to the cyclic sieving phenomenon (CSP) of Reiner--Stanton--White,
which we now recall.
\begin{Definition}{{\cite{MR2087303}}}\label{def:csp1}
  Suppose $C_n$ is a cyclic group of order $n$ generated
  by $\sigma_n$, $W$ is a finite set on which $C_n$ acts,
  and $f(q) \in \bZ_{\geq 0}[q]$. We say the triple
  $(W, C_n, f(q))$ exhibits the
  \textit{cyclic sieving phenomenon (CSP)} if for all
  $r \in \bZ$,
  \begin{align}\label{eq:CSP_eval}
    \begin{split}
      f(\omega_n^r) &= \# W^{\si_n^r} \\
        &\coloneqq \# \{ w \in W: \si_n^r \dd w = w \}
        = \chi^W(\sigma_n^r),
    \end{split}
  \end{align}
  where $ \omega_n $ is a primitive $ n $th root of unity
  and $\chi^W$ is the character of $W$ as a $ C_n $-module.
\end{Definition}
\noindent See \cite{MR2866734} for an excellent
survey and introduction to cyclic sieving. The following
cyclic sieving result also due to Reiner--Stanton--White
is intimately related to \eqref{eq:SYT_evals}. We use $ W^{\stat}(q) $ to denote $ \sum_{ w \in W} q^{\stat(w)} $.
\begin{Theorem}{{\cite[Theorem~8.3, Proposition~4.4]{MR2087303}}}\label{thm:rsw_alpha}
  Let $ \alpha \vDash n $, let $ \W_\al $ denote the set of
  all words of content $ \al $, let $ C_n $ act on
  $ \W_\al $ by rotation, and let
  $ \maj $ denote the major index statistic.
  Then, the triple
  \[
    \lp \W_\alpha, C_n, \W_\alpha^{\maj}(q) \rp
  \]
  exhibits the CSP.
\end{Theorem}
Since the sets $\W_\alpha$ are precisely the $S_n$-orbits for the natural $ S_n $ action on length
$n$ words, \Cref{thm:rsw_alpha} may be thought of as a ``universal
sieving result'' as follows. A very similar observation appeared in
\cite[Prop.~3.1]{MR2837599}.
\begin{Corollary}\label{cor:Sn_univ}
  Let $W$ be a finite set of length $n$ words closed under the
  $S_n$-action. Then, the triple
    \[ \lp W, C_n, W^{\maj}(q) \rp \]
  exhibits the CSP.
\end{Corollary}
In \cite{Ahlbach201837}, the authors introduced a new
statistic on words, $\flex$. As an example,
$\flex(221221) = 2 \cdot 3 = 6$ since $221221$ is the concatenation of
$2$ copies of the primitive word $221$ and $221221$ is third in
lexicographic order amongst its $3$ cyclic rotations. See
\Cref{def:flex} for details. The $\flex$ statistic was designed
to be ``universal'' for cyclic rather than symmetric actions on
words in the following sense.

\begin{Lemma}{{\cite[Lemma~8.3]{Ahlbach201837}}}
  Let $W$ be a finite set of length $n$ words closed under the
  $C_n$-action, where $C_n$ acts by cyclic rotations. Then, the triple
    \[ \lp W, C_n, W^{\flex}(q) \rp \]
  exhibits the CSP.
\end{Lemma}
A corollary of these universal sieving results is
the following equidistribution result. A more refined statement
appeared in \cite{Ahlbach201837}.
\begin{Theorem}\cite[Theorem~8.4]{Ahlbach201837} \label{thm:majn_flex}
  Let $\W_n$ denote the set of length $n$ words, let
  $\maj_n$ denote the major index modulo $n$ taking values
  in $ \{1, \dots, n \} $, and let $\cont$ denote the \textit{content}
  of a word. We then have
    \[ \W_n^{\cont,\maj_n}(\mathbf{x}; q)
       = \W_n^{\cont,\flex}(\mathbf{x}; q). \]
\end{Theorem}
\noindent In \Cref{sec:KW}, we show that the following well-known
result of Kra{\'s}kiewicz--Weyman is essentially a corollary of
\Cref{thm:majn_flex}. Here $\chi^r$ is the linear representation of
the cyclic group $C_n$ given by $ \chi^r(\si_n) = \w_n^r $.

\begin{Theorem}{\cite{MR1867283}}\label{thm:KW}
We have
    \[ \Ch \chi^r\ind_{C_n}^{S_n}
       = \sum_{\lambda \vdash n} a_{\lambda, r} s_\lambda(\mathbf{x}) \]
  where
    \[ a_{\lambda, r} \coloneqq
       \#\{Q \in \SYT(\lambda) : \maj(Q) \equiv_n r\}. \]
\end{Theorem}

Klyachko \cite[Prop.~1]{klyachko74} showed that the Lie modules
$\cL_n$ and the induced representations $\chi^1\ind_{C_n}^{S_n}$
are Schur--Weyl duals. The $\lambda = (n)$ case of Thrall's problem
thus follows from \Cref{thm:KW} when $r=1$. More precisely, Klyachko
expressed both the characteristic of $\chi^1\ind_{C_n}^{S_n}$ and
the character of $\cL_{n}$ as content generating functions on
primitive necklaces of length $n$ words. We generalize this observation
in \Cref{sec:KW} as follows, which also naturally motivates the
$\flex$ statistic.
\begin{Theorem}\label{thm:ind_NFD}
  Let $\NFD_{n, r}$ denote the set of necklaces of length $n$ words
  with \textit{frequency} dividing $r$, $\F_{n, r}$
  denote the set of length $n$ words with $\flex$ equal to $r$, and
  $ \M_{n,r}(\bfx) $ denote the set of length $ n $ words with
  $ \maj_n $ equal to $ r $. Then
    \[ \Ch \chi^r\ind_{C_n}^{S_n} = \NFD_{n, r}^{\cont}(\mathbf{x})
       = \F_{n, r}^{\cont}(\mathbf{x})
        = \M_{n,r}^{\cont}(\bfx). \]
\end{Theorem}
Our new proof of Kra{\'s}kiewicz--Weyman's result reduces the
problem of finding a bijective proof of a well-known
symmetry result following from \Cref{thm:KW} to finding a
bijective proof of the above equidistribution result,
\Cref{thm:majn_flex}; see \Cref{cor:symmetry}. It also provides
a thus far rare example of an instance of cyclic sieving being
used to prove other results rather than vice-versa.

In \Cref{sec:Cn_branching}, we give a new proof of a result of
Stembridge \cite{MR1023791} which settled a conjecture of Stanley
describing the irreducible multiplicities of induced representations
$\chi^r\ind_{\langle\sigma\rangle}^{S_n}$ for arbitrary $ \si \in S_n $.
The corresponding generalized major index statistics arise very naturally
from the combinatorics of orbits and cyclic sieving.

In \Cref{sec:HLM}, we prove and generalize a result of Schocker
\cite{MR1984625} concerning the higher Lie modules. Thrall's
problem may be reduced to the $\lambda = (a^b)$ case by the
Littlewood--Richardson rule. Bergeron--Bergeron--Garsia
\cite{MR1035495} identified the Schur--Weyl dual of
$\cL_{(a^b)}$ as a certain induced module
$\chi^{1, 1}\ind_{C_a \wr S_b}^{S_{ab}}$
where $C_a \wr S_b$ is a wreath product; see
\Cref{ssec:background_wreaths} for details. Schocker gave a
formula for the multiplicity of the irreducible $V^\mu$ in $\cL_{(a^b)}$,
though it involves many subtractions and divisions in general. We
generalize Schocker's formula to all one-dimensional representations
of $C_a \wr S_b$. In our approach, the subtractions and divisions
in Schocker's formula arise naturally from the underlying combinatorics
using M\"obius inversion and Burnside's lemma.

The basic outline of each argument is the same: we obtain an orbit
generating function from an explicit basis of a $\GL(V)$-module, we
construct an appropriate necklace generating function, we use cyclic
sieving to rewrite this generating function using words and descent
statistics like the major index, and we finally apply RSK to get a
Schur expansion. Transitioning from an orbit generating function to
a necklace generating function where we can apply cyclic sieving
involves various combinatorial techniques.

In \Cref{sec:mash}, we discuss applying aspects of
our approach to Thrall's problem in general.
The arguments in the preceding sections strongly suggest attacking
Thrall's problem by considering all branching rules for the inclusion
$C_a \wr S_b \hookrightarrow S_{ab}$ rather than
considering only one such rule. To that end, consider the
irreducible representations $S^{\ul}$ of $C_a \wr S_b$, which are
indexed by the set of $a$-tuples
$ \ul = (\lam^{(1)}, \dots, \lam^{(a)}) $ of partitions with
$ \sum_{r = 1}^{a} |\lam^{(r)}| = b $. We first give the following
plethystic expression for the corresponding characteristic.
\begin{Theorem}\label{thm:casbchars}
  For all integers $ a, b \ge 1$, we have
  \[ \Ch S^{\ul} \ind_{C_a \wr S_b}^{S_{ab}}
     = \prod_{r = 1}^{a}
     s_{\lambda^{(r)}}[\NFD_{a,r}^{\cont}(\mathbf{x})]. \]
\end{Theorem}
 
We then identify the analogues of the $\flex$ and $ \maj_n $ statistics
in this context, which send words to such $ a $-tuples of partitions.
We consequently give the following monomial expansion of the
corresponding graded Frobenius series. See
\Cref{ssec:background_wreaths} and \Cref{sec:mash} for details.
\begin{Theorem}\label{thm:grfrob_flexab}
  Fix integers $ a, b \ge 1 $. We have
  \begin{align*}
    \sum_{\ul}
      \dim S^{\ul} \cdot
      \Ch \lp S^{\ul}\ind_{C_a \wr S_b}^{S_{ab}} \rp q^{\ul}
      &= \W_{ab}^{\cont,\flex_a^b}(\mathbf{x}; q) \\
      &= \W_{ab}^{\cont,\maj_a^b}(\mathbf{x}; q)
  \end{align*}
  where the sum is over all $ a $-tuples $ \ul = (\lam^{(1)}, \dots, \lam^{(a)}) $ of partitions with $ \sum_{r = 1}^{a} |\lam^{(r)}| = b $ and the $q^{\ul}$ are independent indeterminates.
\end{Theorem}

The rest of the paper is organized as follows.
In \Cref{sec:background}, we review combinatorial and
representation-theoretic background. In particular,
we summarize work related to Kra{\'s}kiewicz--Weyman's result,
\Cref{thm:KW}, in \Cref{ssec:background_KW},
and we discuss the current status of Thrall's problem in
\Cref{ssec:background_thrall}. In \Cref{sec:KW},
we present our proof of
Kra{\'s}kiewicz--Weyman's result, \Cref{thm:KW}, using cyclic sieving.
In \Cref{sec:Cn_branching}, we give an analogous proof
of Stembridge's result, \Cref{thm:Stembridge}. In \Cref{sec:HLM}, we give
generalizations of Schocker's result,
\Cref{thm:GeneralizedSchocker}. In
\Cref{sec:mash}, we define the statistics
$\flex_a^b$ and $\maj_a^b$, prove \Cref{thm:casbchars} and
\Cref{thm:grfrob_flexab}, and discuss how the approach
could be used to find the branching rules for
$ C_a \wr S_b \hookrightarrow S_{ab} $.

\section{Background}\label{sec:background}

Here we provide background on
words, tableaux, Schur--Weyl duality, Kra{\'s}kiewicz--Weyman's
result, Thrall's problem, and certain wreath products for use in later
sections. All representations will be over $\bC$.
We write $[n] \coloneqq \{1, \ldots, n\}$, $ \#S $ for the
cardinality of a set $S$, and
\begin{align*}
  \ch{S}{k} &\coloneqq \{\text{all $k$-element subsets of $S$}\}, \\
  \mch{S}{k} &\coloneqq \{\text{all $k$-element multisubsets of $S$}\}.
\end{align*}

\subsection{Words}\label{ssec:background_words}
We now recall standard combinatorial notions
on words and fix some notation. A \textit{word} $w$ of
\textit{length} $ n $ is a sequence
$ w = w_1 w_2 \cdots w_n $ of \textit{letters}
$w_i \in \bZ_{\ge 1}$. The \textit{descent set} of $w$ is
$\Des(w) \coloneqq \{1 \leq i < n : w_i > w_{i+1} \} $.
The \textit{major index} of $w$ is
$\maj(w) \coloneqq \sum_{i \in \Des(w)} i$. Let $ \maj_n(w) $ denote $ \maj(w) $ modulo $ n $ taking values in $ [n] $.

The \textit{content} of a word $w$, written $\cont(w)$,
is the sequence $\alpha = (\alpha_1, \alpha_2, \ldots)$
where $\alpha_j$ is the number of $j$'s in $w$.
Such a sequence $\alpha$ is called a (weak)
\textit{composition} of $n$, written $\alpha \vDash n$. For $ n \ge 1 $ and $ \al \vDash n $, we
write the set of words of length $n$ or content $\alpha$ as
\begin{align*}
  \W_n &\coloneqq \{w = w_1 \cdots w_n : w_i \in \bZ_{\geq 1}\}, \\
  \W_\alpha &\coloneqq \{w \in \W_n : \cont(w) = \alpha \}.
\end{align*}

The set of all words with letters from $\bZ_{\geq 1}$
is a monoid under concatenation. A word is
\textit{primitive} if it is not a power of a smaller
word. Any non-empty word $w$ may be written uniquely
as $w = v^f$ for $f \geq 1$ with $ v $ primitive.
The \textit{period} of $w$,
denoted $ \period(w) $, is the length of $v$. The \textit{frequency}
of $w$, denoted $\freq(w)$, is $f$.

The symmetric group $ S_n $ acts on $ \W_n $ by permuting the letters
according to
\begin{align}\label{eq:rotaction}
  \si \cdot w_1 w_2 \cdots w_n
    \coloneqq w_{\si^{-1}(1)} w_{\si^{-1}(2)} \cdots w_{\si^{-1}(n)}
\end{align}
for all $ \si \in S_n $.
In particular, letting $\sigma_n \coloneqq
(1 \ 2 \ \cdots \ n) \in S_n $ and
$ C_n \coloneqq \langle \si_n \rangle $, the cyclic group
$ C_n $ acts on $ \W_n $ by rotation according to
\begin{equation*}
    \si_n \cdot w_1 w_2 \cdots w_n
    \coloneqq w_{n} w_{1} \cdots w_{n-1}.
\end{equation*}

\begin{Definition}\label{def:necklaces}
  An orbit of $ w \in \W_n $ under rotation is a
  \emph{necklace}, denoted $[w]$. Note that
  $\period(w) = \#[w] $ and $ \freq(w) \dd
  \period(w) = n $. Content, primitivity, period,
  and frequency are all well-defined on necklaces. For
  $ n \ge 1 $, we write
  \begin{align*}
    \N_n   &\coloneqq \{ \tx{necklaces of length $ n $ words} \}.
  \end{align*}
\end{Definition}

\begin{Example}\label{ex:wordandnecklace}
  Consider $ w = 15531553 \in \W_8 $. Then, the length of $w$ is $8$,
  $\Des(w) = \{ 3, 4, 7 \} $, $\maj(w) = 14$,
  and $ \cont(w) = (2,0,2,0,4) $, so $ w \in \W_{(2,0,2,0,4)} $.
  Since $w = 15531553 = (1553)^2 $ and $ 1553 $ is primitive,
  $ w $ is not primitive, $ \period(w) = 4 $, and
  $ \freq(w) = 2$. The necklace of $ w $ is
    \[ [w] = \{15531553, 55315531, 53155315, 31553155 \}
       \in \N_8. \]
\end{Example}

We now recall the $\flex$ statistic from \cite{Ahlbach201837}.
\begin{Definition}\label{def:flex}
  Given $ w \in \W_n $, let $ \lex(w) $ denote the position
  at which $ w $ appears in the lexicographic order of its
  rotations, starting at $1$. The \textit{flex} statistic
  is given by
    \[ \flex(w) = \freq(w) \dd \lex(w). \]
\end{Definition}
\begin{Example}
  If $ w = 21132113 $, its necklace is
  \[ [w] = \{ 11321132, 13211321, 21132113, 32113211 \} \]
  listed in lexicographic order. Since $w$ is in the third position,
  $ \lex(w) = 3 $. Here $ \freq(w) = 2 $, so $ \flex(w) = 6 $.
\end{Example}

\subsection{Generating Functions}\label{ssec:background_gf}

In most triples $ (W, C_n, f(q)) $ that have been found to
exhibit the CSP, $ f(q) $ is a statistic generating function
on $ W $ for some well-known statistic. Given
$ \stat \colon W \to \bZ_{\geq 0} $, we write the corresponding
generating function as
\begin{align*}
  W^{\stat}(q) \coloneqq \sum_{w \in W} q^{\stat(w)}.
\end{align*}
We use natural multivariable analogues of this notation
as well. For example, letting $ \mathbf{x} = ( x_1, x_2, \dots )  $,
\begin{align*}
  \W_n^{\cont, \maj}(\mathbf{x}; q) \coloneqq \sum_{w \in \W_n}
    \mathbf{x}^{\cont(w)} q^{\maj(w)}
    \in \bZ_{\geq 0}[[x_1, x_2, \ldots]][q]
\end{align*}
where $ \mathbf{x}^{(\al_1, \dots, \al_m)}  \coloneqq x_1^{\al_1} \cdots x_m^{\al_m}$.

\subsection{Tableaux}\label{ssec:background_tableaux}

A \textit{partition} of $n$, denoted $\lambda \vdash n$, is
a composition of $n$ whose parts weakly decrease. Write $\Par$
for the set of all partitions. The
\textit{Young diagram} of $\lambda$ is the upper-left justified
collection of cells with $\lambda_i$ entries in the $i$th row
starting from the top. We may write a partition in exponential form
as $\lambda = 1^{m_1} 2^{m_2} \cdots \vdash n$ where $m_i$ is the
number of parts of $\lambda$ of size $i$. In this case, the number
of elements of $S_n$ with cycle type $\lambda$ is $\frac{n!}{z_\lambda}$
where $z_\lambda \coloneqq 1^{m_1} 2^{m_2} \cdots m_1! m_2! \cdots$.

A \textit{semistandard Young tableau}
of shape $\lambda$ is a filling of the Young diagram of $\lambda$
with entries from $\bZ_{\geq 1}$ which weakly increases along
rows and strictly increases along columns. The set of semistandard
Young tableaux of shape $\lambda$ is denoted $\SSYT(\lambda)$.
The \textit{content} of $P \in \SSYT(\lambda)$, denoted $\cont(P)$,
is the composition whose $j$-th entry is the number of $j$'s in $P$.
The set of \textit{standard Young tableaux} of shape $\lambda$,
denoted $\SYT(\lambda)$, is the subset of $\SSYT(\lambda)$
consisting of tableaux of content $(1, \ldots, 1) \vDash n $. The
\textit{descent set} of a tableau $Q \in \SYT(\lambda)$, denoted
$\Des(Q)$, is the set of all $i \in [n-1]$ such that $i+1$
lies in a lower row of $Q$ than $i$.

\begin{Example}
  We draw our tableaux in English notation.
  The semistandard tableau
    \[ P = \Yvcentermath1{\young(112334,23446,3)}
       \in \SSYT(6,5,1) \]
  has $ \cont(P) = (2,2,4,3,0,1) $. The standard tableau
    \[ Q = \Yvcentermath1{\young(125,34,6)}
       \in \SYT(3,2,1) \]
  has $ \Des(Q) = \{ 2,5\} $, and $ \maj(Q) = 7 $.
\end{Example}

Let $ \mathbf{x} = (x_1, x_2, \dots ) $. For a partition $ \lam $, the
\textit{Schur function} $ s_\lam $ is the content generating function
on semistandard tableaux of shape $\lambda$,
  \[ s_{\lam}(\mathbf{x}) \coloneqq \SSYT(\lam)^{\cont}(\mathbf{x})
     \coloneqq \sum_{P \in \SSYT(\lam)} \mathbf{x}^{\cont(P)}. \]
The Schur functions are symmetric in the sense that
they are unchanged under any permutation of the underlying
variables. Two important instances of Schur functions are the
complete homogeneous symmetric functions
\begin{align}\label{eq:hn}
    h_n(\bfx) \coloneqq s_{(n)}(\bfx) = \sum_{i_1 \le \cdots \le i_n} x_{i_1} \cdots x_{i_n}
\end{align}
and the elementary symmetric functions
\begin{align}\label{eq:en}
    e_n(\bfx) \coloneqq s_{(1^n)}(\bfx) = \sum_{i_1 < \cdots < i_n} x_{i_1} \cdots x_{i_n}.
\end{align}
The \textit{power-sum symmetric functions} are given by
  \[ p_n(\bfx)  \coloneqq x_1^n + x_2 ^n + \cdots \quad\text{and}\quad
     p_{(\lam_1, \dots, \lam_k)}(\bfx) \coloneqq p_{\lam_1}(\bfx) \cdots p_{\lam_k}(\bfx). \]

\begin{Definition}
  The Robinson--Schensted--Knuth (RSK) correspondence is a bijection
  \begin{align*}
    \RSK \colon \W_n &\to \bigsqcup_{\lam \vdash n}
      \SSYT(\lam) \x \SYT(\lam), \\
    w &\mapsto (P(w), Q(w)).
  \end{align*}
  The \textit{shape} of $w$ under RSK, denoted $\sh(w)$,
  is the common shape of $P(w)$ and $Q(w)$. Two well-known
  properties of the RSK correspondence are
  \begin{align}\label{eq:RSKprop}
    \cont(w) = \cont(P(w)), \qquad \Des(w) = \Des(Q(w)).
  \end{align}
The fact that $\Des(w) = \Des(Q(w))$ is originally due to
Sch\"utzenberger \cite[Remarque~2]{MR0190017}. See
\cite[Lemma 7.23.1]{MR1676282} for a proof of \eqref{eq:RSKprop}
in the decisive permutation case and \cite[p.404]{MR1676282} for
further historical remarks. See \cite[Chapter 3]{MR1824028}
for more details on RSK.
\end{Definition}

We will repeatedly use the RSK correspondence to transition from the
monomial to the Schur basis. These arguments all rely on the following result.

\begin{Lemma}\label{lem:RSK_Des}
  Suppose $D \subset [n-1]$ and let
    \[ \W_{n, D} \coloneqq \{w \in \W_n : \Des(w) = D\} \]
  be the set of length $n$ words with descent set $D$.
  For $\lambda \vdash n$, let
    \[ a_\lambda^D \coloneqq
       \#\{Q \in \SYT(\lambda) : \Des(Q) = D\}. \]
  Then
    \[ \W_{n, D}^{\cont}(\mathbf{x})
       = \sum_{\lambda \vdash n} a_\lambda^D s_\lambda(\mathbf{x}).
    \]

  \begin{proof}
    Using RSK and \eqref{eq:RSKprop}, we have
    \begin{align*}
      \W_{n, D}^{\cont}(\mathbf{x})
        &= \{w \in \W_n : \Des(w) = D\}^{\cont}(\mathbf{x}) \\
        &= \sum_{\lambda \vdash n}
           \, \sum_{\substack{Q \in \SYT(\lambda)\\\Des(Q) = D}}
           \, \sum_{P \in \SSYT(\lambda)} x^{\cont(P)} \\
        &= \sum_{\lambda \vdash n}
           \, \sum_{\substack{Q \in \SYT(\lambda)\\\Des(Q) = D}}
           s_\lambda(\mathbf{x}) \\
        &= \sum_{\lambda \vdash n} 
           a_\lambda^D s_\lambda(\mathbf{x}).
    \end{align*}
  \end{proof}
\end{Lemma}

\subsection{Schur--Weyl Duality}\label{ssec:background_schurweyl}
We next summarize a few key points from the representation theory
of $S_n$ and $\GL(\bC^m)$. See \cite{MR1464693} for more.

The complex irreducible inequivalent representations
of $S_n$ are canonically indexed by partitions
$\lambda \vdash n$ and are called Specht modules, written
$S^\lambda$. The \textit{Frobenius characteristic map} $\Ch $
is defined by $\Ch S^\lambda \coloneqq s_\lambda(\mathbf{x})$ and
is extended additively to all $S_n$-representations. Since Schur
functions are $\bZ$-linearly independent, computing the irreducible
decomposition of an $S_n$-module $M$ corresponds to computing
the Schur expansion of $\Ch M$.

Let $V$ be a complex vector space of dimension $m$. Endow
$V^{\otimes n}$ with the diagonal left $\GL(V)$-action
and the natural right $S_n$-action given by
permutation of indexes. Given any $S_n$-module $ M $,
define a corresponding $\GL(V)$-module by
\begin{align*}
  E(M) \coloneqq V^{\otimes n} \otimes_{\bC S_n} M,
\end{align*}
which we call the \textit{Schur-Weyl dual} of $M$. The
irreducible inequivalent polynomial representations of
$\GL(V)$ are precisely the Schur Weyl duals of all $ S^\lambda $ where $\lambda$ is a partition with at most $\dim(V)$
non-zero parts \cite[Thm.~8.2.2]{MR1464693}.

Let $E$ be a finite-dimensional, polynomial representation
of $\GL(V)$ and pick a basis $ \{ v_1, \ldots, v_m \} $
for $V$. The \textit{Schur character} of $E$, denoted
$\Ch E$, is the trace of the action of
$ \diag(x_1, \dots, x_m) \in \GL(V) $ on $ E $, where the diagonal
matrix is with respect to the basis $ v_1, \dots, v_m $.
Polynomiality of $E$ implies $ \Ch E \in \bC[x_1, \ldots, x_m]$.
Moreover, $ \Ch(E) $ is a symmetric function of $ x_1, \dots, x_m $.
In fact,
  \[ \Ch V^\lambda = \Ch E(S^\lambda)
     = s_\lambda(x_1, \ldots, x_m, 0, 0, \ldots). \]
Thus, for any $S_n$-module $ M $, we have
  \[ \lim_{m \to \8} \Ch E(M) = \Ch M. \]
In light of this, we often leave dependence on $m$ or
$V$ implicit.

\subsection{Kra{\'s}kiewicz--Weyman Symmetric Functions}\label{ssec:background_KW}
The symmetric functions appearing in \Cref{thm:majn_flex} have
a wealth of important interpretations. Here we summarize some of these
interpretations.

\begin{Definition}
  For $n \geq 1$, let
  \begin{equation}\label{eq:KW.0}
    \KW_n(\mathbf{x}; q) \coloneqq
      \sum_{\substack{\lambda \vdash n\\r \in [n]}}
      a_{\lambda, r} s_\lambda(\mathbf{x}) \, q^r
  \end{equation}
  where $a_{\lambda, r} \coloneqq \#\{Q \in \SYT(\lambda) :
  \maj(Q) \equiv_n r\}$. We call $\KW_n(\bfx; q)$ the $n$th
  \textit{Kra{\'s}kiewicz--Weyman} symmetric function.
\end{Definition}
These symmetric functions are intimately related to the irreducible
representations of certain cyclic groups.
\begin{Definition}\label{def:cyclicirreps}
  Recall $ \si_n \coloneqq (1 \, 2 \, \cdots \, n) \in S_n $ and
  $ C_n \coloneqq \langle \si_n \rangle \le S_n $ be the cyclic
  group of order $ n $ it generates. Fixing any primitive
  $n$th root of unity $\omega_n$, write the irreducible
  characters of $ C_n $ as $ \chi^1, \dots, \chi^n $ where
    \[ \chi^r(\si_n) \coloneqq \omega_n^r. \]
  We sometimes write $\chi^r_n$ if we want to specify the cyclic
  group $ C_n $ as well.
\end{Definition}

\Cref{thm:KW} gives our first interpretation of
$\KW_n(\bfx ; q) $,
\begin{equation}\label{eq:KW.1}
  \KW_n(\mathbf{x}; q)
    = \sum_{r=1}^n \Ch \chi^r\ind_{C_n}^{S_n} q^r.
\end{equation}
Since the regular representation of
$C_n$ is $\oplus_{r=1}^n \chi^r$, when $q=1$ the right-hand side
of \eqref{eq:KW.1} is the Frobenius characteristic of the regular
representation of $S_n$, denoted $ \bC S_n $. The right-hand side of \eqref{eq:KW.1}
is hence similar to a graded Frobenius series for $\bC S_n$ and
tracks branching rules for the inclusion $C_n \hookrightarrow S_n$.
By \Cref{thm:ind_NFD}, we can also write this series as
\begin{align}\label{eq:KWNFD}
  \KW_n(\mathbf{x}; q) = \sum_{r=1}^n \NFD_{n,r}^{\cont}(\bfx) \, q^r.
\end{align}

Now consider the action of $\sigma_n$ on the $S_n$-irreducible
$S^\lambda$. Since $\sigma_n^n = 1 \in S_n $, the action of $\sigma_n$ on
$S^\lambda$ is diagonal with eigenvalues $\omega_n^{k_1},
\omega_n^{k_2}, \ldots$ where $\omega_n$ is a fixed primitive
$n$th root of unity and $1 \leq k_i \leq n$ for each $i$. Let
$P_\lambda(q) \coloneqq q^{k_1} + q^{k_2} + \cdots$ be the
generating function of the \textit{cyclic exponents} $k_1,
k_2, \ldots$, which were studied extensively by Stembridge
\cite{MR1023791}. Using the right-hand side of \eqref{eq:KW.1} and
Frobenius reciprocity quickly gives the following.
\begin{Theorem}[See {\cite[Prop.~1.2, Thm.~3.3]{MR1023791}}]
  The cyclic exponent generating function for $S_n$ is given by
  \begin{equation}\label{eq:KW.1.5}
    \KW_n(\mathbf{x}; q) = \sum_{\lambda \vdash n}
      P_\lambda(q) s_\lambda(\mathbf{x}).
  \end{equation}
\end{Theorem}

Next, extend the regular representation $\bC S_n$ to an
$S_n \times C_n$-module by letting $S_n$ act on the left and $C_n$
act on the right. There is a
straightforward notion of an $S_n \times C_n$-Frobenius
characteristic map given by sending an irreducible
$S^\lambda \otimes \chi^r$ to $s_\lambda(\mathbf{x}) q^r$
where $q$ is an indeterminate satisfying $q^n = 1$.
The following now follows easily using the right-hand side of
\eqref{eq:KW.1}.

\begin{Corollary}{{\cite{MR1867283}}}\label{cor:kw.grfrob}
  The $S_n \times C_n$-Frobenius characteristic of
  the regular representation is
  \begin{equation}\label{eq:KW.2}
    \KW_n(\mathbf{x}; q) = \Ch_{S_n \times C_n} \bC S_n.
  \end{equation}
\end{Corollary}

It is well-known that the type $A_{n-1}$ coinvariant algebra $R_n$
is a graded $S_n$-module which is isomorphic as an
\textit{ungraded} $S_n$-module to $ \bC S_n $.
We may give $R_n$ an $S_n \times C_n$-module structure
by letting $C_n$ act on the $k$th degree component of
$R_n$ by $\sigma_n \cdot f \coloneqq \omega_n^k f$, where
$\omega_n$ is a fixed primitive $n$th root of unity.
Springer and, independently, Kra{\'s}kiewicz--Weyman showed that
$\bC S_n$ and $R_n$ are isomorphic as $S_n \times C_n$-modules.
Consequently, from the right-hand side of \eqref{eq:KW.2},
we have the following.

\begin{Theorem}[Springer {\cite[Prop.~4.5]{MR0354894}}; cf.~ {\cite[Thm.~1]{MR1867283}}]\label{thm:KW_coinvariant}
  The $S_n \times C_n$-Frobenius characteristic of the coinvariant
  algebra $R_n$ is
  \begin{equation}\label{eq:KW.2.5}
    \KW_n(\mathbf{x}; q) = \Ch_{S_n \times C_n} R_n.
  \end{equation}
\end{Theorem}

The graded Frobenius characteristic of the coinvariant algebra
is the modified Hall--Littlewood symmetric function
$\widetilde{Q}_{(1^n)}(\mathbf{x}; q)$ \cite[(I.8)]{MR1168926}.
Consequently, \eqref{eq:KW.2.5} gives
\begin{equation}\label{KW.3}
  \KW_n(\mathbf{x}; q) \equiv \widetilde{Q}_{(1^n)}(\mathbf{x}; q)
  \qquad\text{(mod $q^n-1$)}.
\end{equation}
See also \cite[\S3]{MR2574838} for a nice summary of this connection.

We may instead use the right-hand side of \eqref{eq:KW.0} as a starting
point. From \Cref{lem:RSK_Des}, it follows that
\begin{equation}\label{eq:KW.4}
  \KW_n(\mathbf{x}; q) = \W_n^{\cont, \maj_n}(\mathbf{x}; q).
\end{equation}
\noindent From \Cref{thm:majn_flex} and \eqref{eq:KW.4}, our final
interpretation of $\KW_n(\bfx;q)$ in this subsection is
\begin{equation}\label{eq:KW.5}
  \KW_n(\mathbf{x}; q) = \W_n^{\cont, \flex}(\mathbf{x}; q).
\end{equation}

\subsection{Thrall's Problem}\label{ssec:background_thrall}
We next define the Lie modules $\cL_\lambda$ and summarize
the status of Thrall's problem. See \cite{MR1231799} for more
details.

The \textit{tensor algebra} of $V$ is $T(V)
\coloneqq \oplus_{n=0}^\infty V^{\otimes n}$, which is naturally
a graded $\GL(V)$-representation. Let $\cL(V)$ be the Lie
subalgebra of $T(V)$ generated by $V$, called the
\textit{free Lie algebra} on $V$, so that $\cL(V)$ is a
graded $\GL(V)$-representation with graded components
$\cL_n(V) = V^{\otimes n} \cap \cL(V)$ called \textit{Lie modules}. The
\textit{universal enveloping algebra} $\fU(\cL(V))$ is isomorphic to
$T(V)$ itself. By the Poincar\'e--Birkhoff--Witt Theorem,
\begin{equation*}
  \fU(\cL(V)) \cong \bigoplus_{\lambda
    = 1^{m_1} 2^{m_2} \cdots} \Sym^{m_1}(\cL_1(V))
      \otimes \Sym^{m_2}(\cL_2(V)) \otimes \cdots
\end{equation*}
as graded $ \GL(V) $-representations, where the sum is
over all partitions and $\Sym^m(M)$ is the $ m $th symmetric power  of $ M $ \cite[Lemma 8.22]{MR1231799}. The
\textit{higher Lie module} associated to
$\lambda = 1^{m_1}2^{m_2} \cdots$ is defined to be
\begin{align}\label{eq:llambda}
  \cL_\lambda(V) \coloneqq \Sym^{m_1}(\cL_1(V)) \otimes
    \Sym^{m_2}(\cL_2(V)) \otimes \cdots.
\end{align}
The Lie modules hence yield a $\GL(V)$-module
decomposition $T(V) \cong \oplus_{\lambda \in \Par} \cL_\lambda(V)$.

\textit{Thrall's problem} is the determination of the
multiplicity of $V^\mu$ in $\cL_\lambda(V)$, for instance by
counting explicit combinatorial objects. The well-known
Littlewood--Richardson rule solves the analogous problem for
$V^\mu \otimes V^\nu$. It
follows from \eqref{eq:llambda} and the Littlewood--Richardson
rule that, for the purposes of Thrall's problem, we may
restrict our attention to the case when $\lambda = (a^b)$
is a rectangle. Since
\begin{equation}\label{eq:Lie_rect}
  \cL_{(a^b)}(V) = \Sym^b(\cL_{(a)}(V)),
\end{equation}
the single-row case is particularly fundamental.

Hall \cite[Lemma~11.2.1]{MR0103215} introduced what is now
called the \textit{Hall basis}
for $\cL_{n}(V)$, which, in the $m \to \infty$ limit,
is in content-preserving bijection with primitive necklaces $\NFD_{n,1}$.
For each primitive necklace, Hall associates a bracketing of its
elements using what is now known as the Lyndon factorization
\cite{MR0102539}. He gives an explicit, though computationally
complex, algorithm to express any bracketing as a linear combination
of the bracketings associated to primitive necklaces. Linear
independence of these generators follows from a dimension count.

Klyachko consequently observed that the Schur character of $\cL_{n}$ is
the corresponding content generating function
$\NFD_{n,1}^{\cont}(\mathbf{x})$. Taking symmetric powers, it
follows that in the $m \to \infty$ limit,
$\cL_{(a^b)}(V)$ has a basis indexed by
multisets of primitive necklaces and the Schur
character is the following content generating
function.

\begin{Lemma}[See {\cite[Proposition 1]{klyachko74}}\label{lem:lie_necklaces}]
  We have, in the $ m \to \8 $ limit,
  \begin{align*}
    \Ch \cL_{(a)} = \NFD_{a,1}^{\cont}(\mathbf{x}) \quad \tx{ and } \quad \Ch \cL_{(a^b)} = \mch{\NFD_{a,1}}{b}^{\cont}(\mathbf{x}).
  \end{align*}
\end{Lemma}

\noindent One formulation of Thrall's problem is hence to find
the Schur expansion of the expressions in \Cref{lem:lie_necklaces}.

While we will not have direct need of it, we would be remiss if we did
not mention the following beautiful and important result of Gessel and
Reutenauer \cite[(2.1)]{MR1245159}. The expansion of
$ \Ch \cL_{\lam} $ in terms of Gessel's fundamental quasisymmetric
functions is
\begin{align}\label{eq:GRFexpan}
  \Ch \cL_{\lam}
    = \sum_{ \substack{\si \in S_n \\ \si
        \tx{ has cycle type } \lam } } F_{n, \Des(\si)}(\mathbf{x}),
\end{align}
where
  \[ F_{n,D}(\mathbf{x})
     = \sum_{ \substack{ i_1 \le \dots \le i_n \\
       i_j < i_{j + 1} \tx{ if } j \in D } } x_{i_1} \dots x_{i_n}. \]
Gessel and Reutenauer gave an elegant bijective proof of
\eqref{eq:GRFexpan} in \cite{MR1245159} involving
multisets of primitive necklaces as in \Cref{lem:lie_necklaces}.
Another formulation of Thrall's problem is thus to convert the
right-hand side of \eqref{eq:GRFexpan} to the Schur basis.

Klyachko \cite{klyachko74} was the first to observe the intimate
connections between Lie modules and the linear representations
$\chi^r$ in \Cref{def:cyclicirreps}. Klyachko proved the $r=1$ case
of \Cref{thm:ind_NFD}, that
$\Ch \chi^1\ind_{C_n}^{S_n} = \NFD_{n,1}^{\cont}(\mathbf{x})$.
Combining Klyachko's result, \Cref{lem:lie_necklaces}, the
$r=1$ case of \Cref{thm:ind_NFD}, and Kra{\'s}kiewicz--Weyman's result, \Cref{thm:KW}, solves Thrall's problem when $\lambda = (n)$.
Recall that if $\lambda \vdash n$, then
$a_{\lambda, r} \coloneqq \#\{Q \in \SYT(\lambda) : \maj(Q)
\equiv_n r\}$.

\begin{Corollary}\label{cor:kwthrall}
  For all $\lambda \vdash n \geq 1$, the multiplicity of $V^\lambda$
  in $\cL_{(n)}$ is $a_{\lambda, 1}$.
\end{Corollary}

\noindent Since $\chi^r\ind_{C_n}^{S_n}$ depends up to
isomorphism only on $n$ and $\gcd(n, r)$, we also have the following
well-known symmetry.

\begin{Corollary}\label{cor:alamrsym}
  For all $ \lambda \vdash n \ge 1 $, we have
  $ a_{\lam,r} = a_{\lam, \gcd(n,r)} $.
\end{Corollary}

\begin{Remark} A bijective proof of this symmetry is currently unknown.
\end{Remark}

Thrall's problem is an instance of a plethysm problem as we
next describe. See \cite[Appendix 2]{MR1676282} for more details.
Given polynomial representations of general linear groups
  \[ \rho \colon \GL(V) \to \GL(W) \qquad\text{and}\qquad
     \tau \colon \GL(W) \to \GL(X) \]
where $V, W, X$ are finite-dimensional complex vector spaces,
the \textit{plethysm} of their Schur characters is the Schur
character of their composite:
  \[ (\Ch \tau)[\Ch \rho] \coloneqq \Ch(\tau \circ \rho). \]
It is easy to see that $\Ch \Sym^b(W) = h_b(x_1, \ldots, x_m)$
where $m = \dim(W)$. Consequently, \eqref{eq:Lie_rect} gives
\begin{equation}
  \Ch \cL_{(a^b)} = h_b[\Ch \cL_a].
\end{equation}
Yet another formulation of Thrall's problem is thus
to expand $h_b[\Ch \cL_a]$ in the Schur basis.
Such plethysm problems are notoriously difficult.
However, a combinatorial description for the Schur expansion of
$(\Ch \cL_a)[h_\nu]$ in terms of the \textit{charge} statistic was
given by Lascoux--Leclerc--Thibon in \cite[Thm.~4.2]{MR1261063} and
\cite[Thm.~III.3]{MR1434225}.

\begin{Remark}
At present, Thrall's problem has only been solved in the following
cases:
\begin{itemize}
  \item when $\lambda = (n)$ has a single part
    (see \Cref{cor:kwthrall});
  \item when $\lambda = (1^n)$, $\cL_{(1^n)}$ is the trivial
    representation;
  \item when $\lambda = (2^b)$, $\Ch \cL_{(2^b)} = \sum s_\mu$ where
    the sum is over $\mu \vdash 2b$ with even column sizes
    (see \cite[Ex.~I.8.6(b), p.~138]{MR1354144}).
\end{itemize}
\end{Remark}

\subsection{Wreath Products}\label{ssec:background_wreaths}
The Schur--Weyl duals of the higher Lie modules $\cL_\lambda$ have
also been identified in terms of induced representations of certain
wreath products. Here we summarize this connection as well as some
related aspects of the representation theory of wreath products
which will be used in \Cref{sec:mash}. Our presentation largely
mirrors \cite{MR1023791}.

\begin{Definition}\label{defn:wreathprod}
  Given a group $G$, the \textit{wreath product} of $G$ with $S_n$,
  denoted $G \wr S_n$, is the semidirect product explicitly described
  as follows. $G \wr S_n$ is the set $G^n \times S_n$ with
  multiplication given by
    \[ (g_1, \dots, g_n, \s) \dd (h_1, \dots, h_n, \tau)
       \coloneqq (g_1 h_{\s^{-1}(1)}, \dots,
         g_n h_{\s^{-1}(n)}, \s \tau) \]
  for all $ g_1, \dots g_n, h_1, \dots, h_n \in G $ and
  $ \s, \tau \in S_n $. Furthermore, given $\alpha \vDash n$, set
  $G \wr \prod_i S_{\alpha_i} \coloneqq \prod_i (G \wr S_{\alpha_i})$,
  which has a natural inclusion into $G \wr S_n$. Roughly speaking,
  $G \wr S_n$ can be considered as the group of $n \times n$
  ``pseudo-permutation'' matrices with entries from $G$.

  Now suppose $U$ is a $G$-set and $V$ is an $S_n$-set. There is a
  natural notion of $U \wr V$ as a $G \wr S_n$-set. Explicitly,
  let $U \wr V$ be the set $U^n \times V$ with
  $G \wr S_n$-action given by
    \[ (g_1, \dots, g_n, \si) \dd (u_1, \ldots, u_n, v)
       \coloneqq (g_1 \dd u_{\si^{-1}(1)}, \ldots,
          g_n \dd u_{\si^{-1}(n)}, \si \dd v) \]
  for all $ g_1, \dots, g_n \in G, \si \in S_n, u_1,
  \dots, u_n \in U, v \in V $. There is an analogous notion
  if $U$ is a $G$-module and $V$ is an $S_n$-module,
  namely $U \wr V \coloneqq U^{\otimes n} \otimes V$ with
  $G \wr S_n$-action
    \[ (g_1, \dots, g_b, \si) \dd (u_1 \otimes \dots \otimes u_b \otimes v)
       \coloneqq (g_1 \dd u_{\si^{-1}(1)}) \otimes \dots \otimes
          (g_b \dd u_{\si^{-1}(b)}) \otimes (\si \dd v) \]
  extended $\bC$-linearly.
\end{Definition}

Since $S_a$ acts naturally and faithfully on $[a]$, $[a] \wr 1_b$ has a
natural $S_a \wr S_b$-action, where $1_b$ denotes the trivial
$S_b$-set. Identifying $[a] \wr 1_b$ with the
set $[ab]$ and noting that the action remains faithful gives an
inclusion $S_a \wr S_b \hookrightarrow S_{ab}$. Similarly we have
an inclusion $C_a \wr S_b \hookrightarrow S_{ab}$. More concretely,
$C_a \wr S_b$ acts faithfully on $[ab]$ by permuting the $b$
size-$a$ intervals in $[ab]$ amongst themselves and cyclically
rotating each size-$a$ interval independently.

\begin{Remark}
  The induction product of two symmetric group representations
  corresponds to the product of their Frobenius characteristics,
  so that if $U$ is an $S_a$-module and $V$ is an $S_b$-module, then
  \cite[Prop.~7.18.2]{MR1676282},
  \begin{align}\label{eq:prodchar}
    \Ch \lp U \otimes V\ind_{S_a \times S_b}^{S_{a+b}} \rp
    = (\Ch U)(\Ch V).
  \end{align}
  In \Cref{ssec:background_thrall}, we considered the plethysm of
  Schur characters of general linear group representations.
  The corresponding operation for Frobenius characters of
  symmetric group representations is less well-known and involves
  wreath products as follows.
  Given two symmetric functions $f$ and $g = m_1 + m_2 + \cdots$
  where the $m_i$ are all monomials,
  their \textit{plethysm} is given by \cite[Def.~A2.6]{MR1676282}
  \begin{align}\label{eq:plethysmdefn}
    f[g] \coloneqq f(m_1, m_2, \dots),
  \end{align}
  which is well-defined since $f$ is symmetric. Then, if
  $U$ is an $S_a$-module and $V$ is an $S_b$-module, we have
  (see \cite[Thm.~A2.8]{MR1676282} or
  \cite[Appendix A, (6.2)]{MR1354144})
  \begin{align}\label{eq:plethsymchar}
    \Ch \lp (U \wr V)\ind_{S_a \wr S_b}^{S_{ab}} \rp
      = \Ch(V)[\Ch(U)].
  \end{align}
\end{Remark}

When $G$ is a finite group, Specht \cite{Specht1932Eine}
described the complex inequivalent irreducible representations
of $G \wr S_n$ in terms of those for $G$ and $S_n$,
the conjugacy classes of $G$, and wreath products. In the case
$C_a \wr S_b$, they are indexed by the following objects.

\begin{Theorem}[{\cite{Specht1932Eine}}; see {\cite[Thm.~4.1]{MR1023791}}]\label{thm:CaSb_irreps}
  The complex inequivalent irreducible representations of
  $C_a \wr S_b$ are indexed by $ a $-tuples
  $ \ul =  (\lam^{(1)}, \dots, \lam^{(a)}) $ of partitions
  with $ \sum_{r = 1}^a |\lam^{(r)}| = b $. In particular, they are
  given by
  \begin{equation}\label{eq:specht}
    S^{\ul}
      \coloneqq \lp (\chi_a^1 \wr S^{\lambda^{(1)}}) \otimes
         \dots \otimes (\chi_a^{a} \wr S^{\lambda^{(a)}}) \rp
         \ind_{C_a \wr S_{ \al(\ul) } }^{C_a \wr S_b},
  \end{equation}
  where
  \begin{align*}
    \al(\ul)
      &\coloneqq (|\lambda^{(1)}|, \dots, |\lambda^{(a)}|) \vDash b, \\
    S_{\al(\ul)}
      &\coloneqq S_{|\lambda^{(1)}|} \times \cdots \times
        S_{|\lambda^{(a)}|},
  \end{align*}
$ \chi_a^r $ is as defined in \Cref{def:cyclicirreps}, and $C_a \wr S_{\al(\ul)}$ is viewed naturally as a subgroup
  of $C_a \wr S_b$.
\end{Theorem}

One consequence of \Cref{thm:CaSb_irreps} is
\begin{align}\label{eq:dimcasbirreps}
  \dim(S^{\ul})
    = \ch{b}{\al(\ul)} \prod_{r = 1}^{a} \# \SYT(\lam^{(r)}).
\end{align}
Another consequence is an explicit description of the one-dimensional
representations of $C_a \wr S_b$, which are as follows.
\begin{Definition}
  Fix integers $a, b \geq 1$. Let
    \[ \chi^{r, 1} \coloneqq \chi^r_a \wr 1_b \qquad \text{ and }
       \chi^{r, \epsilon} \coloneqq \chi^r_a \wr \epsilon_b \]
  where $r=1, \ldots, a$ and $1_b$ and $\epsilon_b$ are the trivial
  and sign representations of $S_b$, respectively. When
  $b=1$, $ \epsilon_b = 1_b $, in which case $ \chi^{r, 1}
  = \chi^{r, \epsilon} = \chi_a^r $. We sometimes write
  $\chi^{r, 1}_{(a^b)}$ or $\chi^{r, \epsilon}_{(a^b)}$ if we
  want to specify the group $C_a \wr S_b$ as well.
\end{Definition}

Bergeron--Bergeron--Garsia \cite{MR1035495} extended Klyachko's
observation by showing that the Schur--Weyl dual of $\cL_{(a^b)}$ is
$\chi^{1, 1}\ind_{C_a \wr S_b}^{S_{ab}}$. We next give a different
argument of this fact which is straightforward given the preceding
background and which uses a lemma we will require
later in \Cref{sec:mash}.

\begin{Lemma}\label{lem:higher_schur_weyl}
  We have
    \[
\Ch \chi^{1, 1}\ind_{C_a \wr S_b}^{S_{ab}} = \mch{\NFD_{a,1}}{b}^{\cont}(\mathbf{x}).
\]

  \begin{proof} By \Cref{lem:indwreath} below and the fact that $ C_a \wr S_b \sube S_a \wr S_b \sube S_{ab} $, we have
    \begin{equation}
      \chi^{1, 1}\ind_{C_a \wr S_b}^{S_{ab}}
        = (\chi^1_a \wr 1_b)\ind_{C_a \wr S_b}^{S_{ab}}
          \cong (\chi^1_a\ind_{C_a}^{S_a} \wr 1_b)
          \ind_{S_a \wr S_b}^{S_{ab}}.
    \end{equation}
    By \eqref{eq:plethsymchar} and the $r=1$ case of
    \Cref{thm:ind_NFD},
    \begin{align}
      \Ch (\chi^1_a\ind_{C_a}^{S_a} \wr 1_b)
           \ind_{S_a \wr S_b}^{S_{ab}} = (\Ch 1_b)[\Ch \chi^1_a\ind_{C_a}^{S_a}] = h_b[\NFD_{a,1}^{\cont}(\mathbf{x})],
    \end{align}
 since $ \Ch(1_b) = h_b(\bfx) $. Now, $h_b[\NFD_{a,1}^{\cont}(\mathbf{x})]
    = \mch{\NFD_{a,1}}{b}^{\cont}(\mathbf{x})$ from
    the definition of plethysm, \eqref{eq:plethysmdefn}, and the definition of $ h_b$, \eqref{eq:hn}. The
    result will be complete once we prove \Cref{lem:indwreath}.
  \end{proof}
\end{Lemma}

\begin{Corollary}[{\cite[\S4.4]{MR1035495}}; see also {\cite[Thm.~8.24]{MR1231799}}]\label{cor:higher_schur_weyl}
  The Schur--Weyl dual of $\cL_{(a^b)}$ is
  $\chi_a^1\ind_{C_a \wr S_b}^{S_{ab}}$.

  \begin{proof}
    Combine \Cref{lem:lie_necklaces} and \Cref{lem:higher_schur_weyl}.
  \end{proof}
\end{Corollary}

Indeed, the Schur--Weyl duals of general $\cL_\lambda$ can be expressed
very explicitly in terms of induced linear representations as follows.
Suppose $\sigma \in S_n$ has cycle type $\lambda$. Write $Z_\lambda$
for the centralizer of $\sigma$ in $S_n$. When $\lambda = (a^b)$, it is
straightforward to see that $Z_{(a^b)} \cong C_a \wr S_b$.
Furthermore, when $\lambda = 1^{b_1} 2^{b_2} \cdots k^{b_k}$
is written in exponential notation, we have $Z_\lambda \cong
Z_{(1^{b_1})} \times Z_{(2^{b_2})} \times \cdots \times Z_{(k^{b_k})}$.

\begin{Corollary}[see {\cite[Thm.~8.24]{MR1231799}}]
  Suppose $\lambda = 1^{b_1} 2^{b_2} \cdots k^{b_k} \vdash n$. Let
  $\chi^{1, 1}_\lambda$ denote the linear representation of
  $Z_\lambda \leq S_n$ given by the (outer) tensor product of
  the representations $\chi_{(i^{b_i})}^{1, 1}$ of $C_i \wr S_{b_i}$
  for $1 \leq i \leq k$. Then, the Schur--Weyl dual of $\cL_\lambda$ is
  $\chi^{1, 1}_\lambda\ind_{Z_\lambda}^{S_n}$.

  \begin{proof}
    Using in order \eqref{eq:llambda}, multiplicativity of
    Schur characters under tensor products, \Cref{cor:higher_schur_weyl},
    \eqref{eq:prodchar}, \Cref{lem:indtensor} and transitivity
    of induction, the fact that
    $Z_\lambda \cong \prod_{i=1}^k Z_{(i^{b_i})}$, and the definition of
    $\chi_\lambda^{1, 1}$, we have
    \begin{align*}
      \Ch \cL_\lambda
        &= \Ch\left(\bigotimes_{i=1}^k \cL_{(i^{b_i})}\right)
         = \prod_{i=1}^k \Ch \cL_{(i^{b_i})}
         = \prod_{i=1}^k
           \Ch \chi_{(i^{b_i})}^{1, 1}\ind_{Z_{(i^{b_i})}}^{S_{ib_i}} \\
        &= \Ch
           \left(\bigotimes_{i=1}^k
           \chi_{(i^{b_i})}^{1, 1}\ind_{Z_{(i^{b_i})}}^{S_{ib_i}}\right)
           \ind_{S_{1b_1} \times S_{2b_2} \times \cdots}^{S_n}
         = \Ch\left(\bigotimes_{i=1}^k
           \chi_{(i^{b_i})}^{1, 1}\right)\ind_{Z_\lambda}^{S_n} \\
        &= \Ch \chi_\lambda^{1, 1}\ind_{Z_\lambda}^{S_n}
    \end{align*}
    The result will be complete once we prove \Cref{lem:indtensor}.
  \end{proof}
\end{Corollary}

\begin{Lemma}\label{lem:indwreath}
  Suppose that $H$ is a subgroup of a group $G$,
  that $U$ is an $H$-module, and that $V$ is an $S_n$-module.
  Then
    \[ (U \wr V) \ind_{H \wr S_n}^{G \wr S_n}
       \cong \lp U\ind_H^G\rp\wr V \]
  as $G \wr S_n$-modules.

  \begin{proof}
    As sets, we have
    \begin{align*}
      (U \wr V)\ind_{H \wr S_n}^{G \wr S_n}
        &= \bC(G \wr S_n) \otimes_{\bC(H \wr S_n)}
           (U^{\otimes n} \otimes V), \\
      (U\ind_H^G) \wr V
        &= (\bC G \otimes_{\bC H} U)^{\otimes n} \otimes V.
    \end{align*}
    Define
    \begin{align*}
      \phi \colon (U \wr V)\ind_{H \wr S_n}^{G \wr S_n}
        &\to (U\ind_H^G) \wr V, \\
      \psi \colon (U\ind_H^G) \wr V
        &\to (U \wr V)\ind_{H \wr S_n}^{G \wr S_n}
    \end{align*}
    by
    \begin{align*}
      \phi((g_1, &\ldots, g_n, \tau) \otimes
           (u_1 \otimes \cdots \otimes u_n \otimes v)) \\
        &\coloneqq (g_1 \otimes u_{\tau^{-1}(1)}) \otimes \cdots \otimes
            (g_n \otimes u_{\tau^{-1}(n)}) \otimes (\tau \cdot v), \\
      \psi((g_1 &\otimes x_1) \otimes \cdots
           \otimes (g_n \otimes x_n) \otimes y) \\
        &\coloneqq (g_1, \ldots, g_n, 1) \otimes
            (x_1 \otimes \cdots \otimes x_n \otimes y)
    \end{align*}
    extended $\bC$-linearly. It is straightforward to check directly
    that $\phi$ and $\psi$ are well-defined,
    $G \wr S_n$-equivariant, and mutual inverses. Note that showing $ \psi \circ \phi(x) = x $ requires using the relation
\[
	(g_1, \dots, g_n, \tau) \otimes z = (g_1, \dots, g_n,1) \otimes (\tau \dd z)
\]
in $ \bC(G \wr S_n) \otimes_{\bC(H \wr S_n)} (U^{\otimes n} \otimes V) $ for $ g_1, \dots, g_n \in G, \tau \in S_n, z \in U^{\otimes n} \otimes V $.
\end{proof}
\end{Lemma}

\begin{Lemma}\label{lem:indtensor}
  Suppose that $H_1, \ldots, H_k$ are subgroups of groups
  $G_1, \ldots, G_k$ and that $U_i$ is an $H_i$-module for
  $1 \leq i \leq k$. Then
    \[ \left(U_1 \otimes \cdots \otimes U_k\right)
       \ind_{H_1 \times \cdots \times H_k}^{G_1 \times \cdots \times G_k}
       \cong
       U_1\ind_{H_1}^{G_1} \otimes \cdots \otimes U_k\ind_{H_k}^{G_k}
     \]
  as $G_1 \times \cdots \times G_k$-modules.

  \begin{proof}
    Having chosen bases for both sides, there is a natural $\bC$-linear
    map between them. It is easy to check this is also $G_1 \times \cdots
    \times G_k$-equivariant. The details are omitted.
  \end{proof}
\end{Lemma}

\section{Cyclic Sieving and Kra{\'s}kiewicz--Weyman's Result}\label{sec:KW}
In this section, we first build on work of Klyachko to prove
\Cref{thm:ind_NFD}. We then recover Kra{\'s}kiewicz--Weyman's
result, \Cref{thm:KW}, and discuss some benefits of our approach.

Klyachko observed in \cite[Prop.~1]{klyachko74} that
$E(\chi^1\ind_{C_n}^{S_n})$, like
$\cL_{(n)}$, also has a basis indexed by primitive necklaces.
Klyachko's argument may be readily generalized to
$E(\chi^r\ind_{C_n}^{S_n})$ as follows. Recall from the introduction that
\begin{align*}
  \NFD_{n, r}
    &\coloneqq \{ N \in \N_n : \freq(N) \mid r \}, \\
  \F_{n, r}
    &\coloneqq \{ w \in \W_n : \flex(w) = r \}, \\
  \M_{n,r}
    &\coloneqq \{ w \in \W_n : \maj_n(w) = r \}.
\end{align*}
In particular, $\NFD_{n,n} = \N_n $, and $ \NFD_{n,1} $ is the set of
primitive necklaces of length $ n $.

\begin{Theorem}\label{thm:basis_onerow}
  There is a basis for $E(\chi^r\ind_{C_n}^{S_n})$
  indexed by necklaces of length $n$ words with letters from
  $[m]$ and with frequency dividing $r$. Moreover,
  \begin{align}\label{eq:IndcharNFD}
    \Ch \chi^r\ind_{C_n}^{S_n} = \NFD_{n,r}^{\cont}(\mathbf{x}).
  \end{align}

  \begin{proof}
    Suppose the underlying vector space $ V $ has basis
    $ \{ v_1, \dots, v_m \} $. By a slight abuse of
    notation, we may view $\chi^r$ as the vector space
    $\bC$ with the left $C_n$-action $\sigma_n \cdot 1 \coloneqq
    \omega_n^r$. Since $\chi^r\ind_{C_n}^{S_n} \coloneqq
    \bC S_n \otimes_{\bC C_n} \chi^r$, we have
    \begin{align}
      E(\chi^r\ind_{C_n}^{S_n})
        &= V^{\otimes n} \otimes_{\bC S_n} \bC S_n
           \otimes_{\bC C_n} \chi^r
           \cong V^{\otimes n} \otimes_{\bC C_n} \chi^r
    \end{align}
    where $C_n$ acts on $V^{\otimes n}$ on the right by
    ``rotating'' the components of simple tensors. A spanning
    set for $V^{\otimes n} \otimes_{\bC C_n} \chi^r$ is
    given by all $v_{i_1} \otimes \cdots \otimes v_{i_n} \otimes 1$,
    which we abbreviate as $[i_1\ \cdots\ i_n]$. Acting by
    $\sigma_n^{-1} $ on $\chi^r$ on the left or on $V^{\otimes n}$ on the
    right gives the relation
  \begin{align}
    [i_1\ \cdots\ i_n] = \omega_n^{r} [i_2 \, \cdots \, i_n \, i_{1}].
  \end{align}
    This relation shows that $[i_1\cdots i_n]$ is well-defined
    on the level of necklaces, at least up to nonzero scalar
    multiplication, which explains our notation.
    If the word $i_1\cdots i_n$ has frequency $f$ and period $p$,
    we then find
    \begin{align*}
      [i_1\ \cdots\ i_n]
        &= \frac{1}{n} \sum_{j=0}^{n-1}
           \omega_n^{jr} [i_{j+1}\ \cdots\ i_n\ i_1\ \cdots\ i_j] \\
        &= \frac{1}{n} \sum_{k=0}^{p-1}
           \left(\sum_{\ell=0}^{f-1} \omega_n^{(\ell p + k)r}\right)
           [i_{k+1} \cdots i_n i_1 \cdots i_k] \\
        &= \frac{1}{n} \left(\sum_{\ell=0}^{f-1}
           \omega_n^{\ell pr}\right)
           \sum_{k=0}^{p-1} \omega_n^{kr}
           [i_{k+1} \cdots i_n i_1 \cdots i_k].
    \end{align*}
    Since $\omega_n^p$ is a primitive $n/p=f$-th root of unity,
    the factor $\sum_{\ell=0}^{f-1} \omega_n^{\ell pr}$ is
    nonzero if and only if $\omega_n^{pr} = 1$, so if and only if
    $f \mid r$. Picking representatives for necklaces with frequency
    dividing $r$ thus gives a spanning set for
    $E(\chi^r\ind_{C_n}^{S_n})$, and it is easy to see it is in fact
    a basis. Diagonal matrices act on this basis via
    \begin{align}
      \diag(x_1, \ldots, x_n) \cdot [i_1\ \cdots\ i_n]
        = \bfx^{\cont(i_1\cdots i_n)} [i_1\ \cdots\ i_n],
    \end{align}
    from which it follows that the Schur character is the content
    generating function of necklaces of length $n$ words with letters
    from $[m]$ and with frequency dividing $r$. Letting $m \to \infty$,
    \eqref{eq:IndcharNFD} follows.
  \end{proof}
\end{Theorem}

\begin{Lemma}\label{lem:NFD_F}
  We have
    \[ \NFD_{n, r}^{\cont}(\mathbf{x}) = \F_{n, r}^{\cont}(\mathbf{x}) = \M_{n,r}^{\cont}(\bfx).
    \]

  \begin{proof}
    Consider the map
    \begin{align*}
      \iota \colon \F_{n, r} &\to \NFD_{n,r} \\
      \iota(w) &\coloneqq [w].
    \end{align*}
    Since $\flex(w) = \freq(w) \lex(w) = r$,
    we have $ \freq(w) \mid r$, so $[w] \in \NFD_{n, r}$.
    Thus, $\iota$ is in fact a map from $ \F_{n, r} $ to $ \NFD_{n,r} $. Since each necklace in $ \NFD_{n,r} $ contains exactly one word with $ \flex $ equal to $ r $, $\iota$ is a content-preserving bijection. Therefore, $ \NFD_{n, r}^{\cont}(\bfx) = \F_{n, r}^{\cont}(\bfx) $.

Using \Cref{thm:majn_flex}, we have
\[
	\W_n^{\cont, \flex}(\mathbf{x}; q) = \W_n^{\cont, \maj_n}(\mathbf{x}; q),	
\]
which means $ \F_{n, r}^{\cont}(\mathbf{x}) = \M_{n,r}^{\cont}(\bfx) $.
  \end{proof}
\end{Lemma}

\begin{Remark}
  From \Cref{thm:basis_onerow} and \Cref{lem:NFD_F}, the
  Schur character of $\chi^r\ind_{C_n}^{S_n}$ may be described
  as a content generating function for certain necklaces or
  for certain words. This proves \Cref{thm:ind_NFD} from the
  introduction.
\end{Remark}

We may now present our remarkably direct proof of
Kra{\'s}kiewicz--Weyman's result, \Cref{thm:KW}, using cyclic sieving.

\begin{proof}[Proof (of \Cref{thm:KW}).]
  The argument in \Cref{thm:basis_onerow} exhibited an explicit basis
  of the Schur module $E(\chi^r\ind_{C_n}^{S_n})$, showing that
    \[ \sum_{r=1}^n \Ch \chi^r\ind_{C_n}^{S_n} q^r
       = \sum_{r=1}^n \NFD_{n, r}^{\cont}(\mathbf{x}) \, q^r. \]
  From \Cref{lem:NFD_F}, the bijection $\iota \colon \F_{n, r} \too{\sim}
  \NFD_{n, r}$ given by $w \mapsto [w]$ gives
    \[ \sum_{r=1}^n \NFD_{n, r}^{\cont}(\mathbf{x}) \, q^r
       = \W_n^{\cont, \flex}(\mathbf{x}; q). \]
  Using universal cyclic sieving on words for $S_n$-orbits and
  $C_n$-orbits as described in the introduction,
  \Cref{thm:majn_flex} now gives
    \[ \W_n^{\cont, \flex}(\mathbf{x}; q) =
       \W_n^{\cont, \maj_n}(\mathbf{x}; q). \]
  Using the RSK algorithm, \Cref{lem:RSK_Des} gives
    \[ \W_n^{\cont, \maj_n}(\mathbf{x}; q)
       = \sum_{\substack{\lambda \vdash n \\r \in [n]}}
         a_{\lambda, r} s_\lambda(\mathbf{x}) \, q^r. \]
  Combining all of these equalities and extracting the
  coefficient of $q^r$ gives the result.
\end{proof}

Every step of the preceding proof uses an explicit bijection with
the exception of the appeal to cyclic sieving through
\Cref{thm:majn_flex}. This suggests the problem of finding a
bijective proof of \Cref{thm:majn_flex}.

\begin{Problem}\label{prob:maj_flex}
  For each $n \geq 1$, find an explicit, content-preserving
  bijection
    \[ \phi \colon \W_n \to \W_n \]
  such that $\maj_n(w) = \flex(\phi(w))$.
\end{Problem}

\begin{Corollary}\label{cor:symmetry}
  A solution to \Cref{prob:maj_flex} would yield an explicit, bijective
  proof of the identity
  \begin{equation}\label{eq:maj_flex_sym}
    \sum_{\lambda \vdash n} a_{\lambda, r} s_\lambda(\mathbf{x})
      = \sum_{\lambda \vdash n} a_{\lambda, s}
         s_\lambda(\mathbf{x})
  \end{equation}
  for any $r, s \in \bZ$ where $\gcd(n, r) = \gcd(n, s)$.
  \begin{proof}
    We have content-preserving bijections
      \[ \bigsqcup_{\lambda \vdash n} \SSYT(\lambda) \times
         \{Q \in \SYT(\lambda) : \maj(Q) \equiv_n r\}
         \too{\mathrm{RSK}} \M_{n, r}
         \too{\phi} \F_{n, r}
         \too{\iota} \NFD_{n, r}. \]
    Now note that
      \[ \NFD_{n, r} = \NFD_{n, \gcd(n, r)}
         = \NFD_{n, \gcd(n, s)} = \NFD_{n, s}. \]
    We thus have an explicit, content-preserving bijection
    \begin{align*}
      \bigsqcup_{\lambda \vdash n} \SSYT(\lambda) &\times
         \{Q \in \SYT(\lambda) : \maj(Q) \equiv_n r\} \\
         &\too{\sim} \bigsqcup_{\lambda \vdash n} \SSYT(\lambda) \times
         \{Q \in \SYT(\lambda) : \maj(Q) \equiv_n s\}
    \end{align*}
    from which \eqref{eq:maj_flex_sym} follows.
  \end{proof}
\end{Corollary}

\begin{Remark}
  The most difficult step in our proof of \Cref{thm:KW} is the
  universal $S_n$-cyclic sieving result, \Cref{cor:Sn_univ},
  or equivalently \Cref{thm:rsw_alpha}. The proof in
  \cite{MR2087303} of \Cref{thm:rsw_alpha} perhaps unsurprisingly
  uses several of the interpretations of the Kra{\'s}kiewicz--Weyman
  symmetric functions from \Cref{ssec:background_KW}, in particular
  \Cref{thm:KW_coinvariant} involving $\Ch_{S_n \times C_n} R_n$.
  However, both Kra{\'s}kiewicz--Weyman's
  and Springer's original proofs of \Cref{thm:KW_coinvariant}
  hinge upon \eqref{eq:SYT_evals}. Indeed, Kra{\'s}kiewicz--Weyman
  showed explicitly in \cite[Prop.~3]{MR1867283} that
  $\Ch_{S_n \times C_n} \bC S_n = \Ch_{S_n \times C_n} R_n$
  is easily equivalent to \eqref{eq:SYT_evals}. Springer's argument
  proving \eqref{eq:SYT_evals} uses a Molien-style formula, while
  Kra{\'s}kiewicz--Weyman's argument uses a recursion
  involving $\ell$-cores and skew hooks.

  One may thus ask about the relationship between \eqref{eq:SYT_evals}
  and the cyclic sieving result, \Cref{thm:rsw_alpha}.
  Using stable principal specializations, one can consider earlier
  approaches to have been ``in the $s$-basis'' and our approach to have
  been ``in the $h$-basis'' in the following sense. Let
  $\tau^\lambda$ be the $S_n$-character of $1\ind_{S_\lambda}^{S_n}$,
  which has $ \Ch(1\ind_{S_\lambda}^{S_n}) = h_\lam $. We have
  \begin{alignat*}{3}
    \chi^\lambda(\sigma_n^r)
      &= \SYT(\lambda)^{\maj}(\omega_n^r)
      &&= (1-q)\cdots(1-q^n)
         s_\lambda(1, q, q^2, \ldots)|_{q=\omega_n^r}, \\
    \tau^\lambda(\sigma_n^r)
      &= \W_\lambda^{\maj}(\omega_n^r)
      &&= (1-q)\cdots(1-q^n)
         h_\lambda(1, q, q^2, \ldots)|_{q=\omega_n^r}.
  \end{alignat*}
  where the first equality is \eqref{eq:SYT_evals}, the
  second is \cite[Prop.~7.19.11]{MR1676282}, the third is
  \Cref{thm:rsw_alpha}, and the fourth is \cite[Prop.~7.8.3]{MR1676282}
  and \cite[Art.~6]{MR1506186}. Our approach suggests that, as far
  as the Kra{\'s}kiewicz--Weyman theorem is concerned, the $h$-basis
  arises more directly.

  In \cite{Ahlbach201837}, the authors proved a refinement of
  \Cref{thm:rsw_alpha}. Since earlier approaches to
  \Cref{thm:rsw_alpha} involving representation theory could not
  readily be adapted to this refinement, the argument instead
  uses completely different and highly combinatorial techniques.
  Thus, the arguments in \cite{Ahlbach201837} and the proof
  of \Cref{thm:KW} together give an essentially self-contained
  proof of Kra{\'s}kiewicz--Weyman's result.
\end{Remark}

\section{Induced Representations of Arbitrary Cyclic Subgroups of \texorpdfstring{$S_n$}{Sn}}\label{sec:Cn_branching}
We next generalize the discussion in \Cref{sec:KW} to branching
rules for general inclusions
$\langle \sigma\rangle \hookrightarrow S_n$, recovering a result
of Stembridge, \Cref{thm:Stembridge}. Following the outline
of the previous section, we express the relevant
characters in turn as a certain orbit generating
function, \Cref{thm:basis_cyclic}, a
necklace generating function, \Cref{lem:OIRNFexpan},
and a generating function on words, \Cref{lem:OIRNFDexpan}.
Two variations on the major index, $\bfmaj_\nu$ and $\maj_\nu$,
arise quite naturally from our argument. The CSP
\Cref{thm:rsw_alpha} again plays a decisive role.

Throughout this section, let $\sigma \in S_n$, let $C$ be the
cyclic group generated by $\sigma$, and let $ \ell \coloneqq \# C $
be the order of $\sigma$. Fixing a primitive $\ell$-th root of
unity $\omega_\ell$, let $ \chi^r \colon C \to \bC $
for $r=1, \ldots, \ell$ be the linear $C$-module given by
$ \chi^r(\sigma) \coloneqq \omega_\ell^r $. We begin by updating our
notation for this setting and generalizing \Cref{thm:basis_onerow}.

\begin{Definition}
  In analogy with \Cref{def:necklaces}, suppose $\cO$ is
  an orbit of $\W_n$ under the restricted $C$-action.
  The \textit{period} of $\cO$ is $\#\cO$
  and the \textit{frequency} of $\cO$, written $\freq(\cO)$, is
  the stabilizer-order of any element of $\cO$, or
  equivalently $\freq(\cO) = \f{\ell}{\#\cO} $.
  The set of orbits of words whose frequency divides $r$ is
    \[ \OFD_{C, r} \coloneqq \{\text{$C$-orbits $\cO$ of $\W_n$}
       : \freq(\cO) \mid r\}. \]
\end{Definition}

\begin{Theorem}\label{thm:basis_cyclic}
  There is a basis for $E(\chi^r\ind_C^{S_n})$
  indexed by $ C $-orbits of length $n$ words with
  letters from $[m]$ and with frequency dividing $r$.
  Moreover,
  \begin{align}
    \Ch \lp \chi^r\ind_C^{S_n} \rp
    = \OFD_{C,r}^{\cont}(\mathbf{x}).
  \end{align}

  \begin{proof}
    The proof of \Cref{thm:basis_onerow} goes through
    verbatim with the $C$-action replacing the $C_n$-action.
  \end{proof}
\end{Theorem}

Our goal is broadly to replace $\OFD_{C, r}^{\cont}(\mathbf{x})$ with
a necklace generating function, apply cyclic sieving
to get a major index generating function on words, and then
apply RSK to get a Schur expansion.

\begin{Notation}
  For the rest of the section, suppose that $\sigma$ has
  disjoint cycle decomposition
  $\sigma = \sigma_1 \cdots \sigma_k$ with
  $\nu_i \coloneqq |\sigma_i|$. Consequently, $\ell
  = |\langle\sigma\rangle| = \lcm(\nu_1, \ldots, \nu_k)$. Further, write
    \[ C_\nu \coloneqq \{\sigma_1^{r_1} \cdots \sigma_k^{r_k} \in S_n
                  : r_1, \ldots, r_k \in \bZ\}
         \cong C_{\nu_1} \times \cdots \times C_{\nu_k} \]
  where $C_{\nu_i} \coloneqq \langle \sigma_i\rangle \subset S_n$.
  Thus, we have $C \subset C_\nu \subset S_n$.
\end{Notation}

In \Cref{sec:KW}, we considered the $C_n$-orbits of $\W_n$, namely
necklaces $N \in \N_n$. The frequency of $N$ is the stabilizer-order
of $N$, i.e.~$\freq(N) = \#\Stab_{C_n}(N)$. We may group together
$C_n$-orbits of $\W_n$ according to their stabilizer sizes by letting
\begin{align}\label{eq:NFdef}
  \NF_{n, r} \coloneqq \{N \in \N_n : \freq(N) = r\}
\end{align}
be the set of necklaces of length $n$ words with frequency $r$.
Similarly, $\NFD_{n, r}$ consists of $C_n$-orbits of $\W_n$
whose stabilizer is contained in the common stabilizer of
$\NF_{n, r}$.

Analogously, the $C_\nu$-orbits of $\W_n$ can be identified with products
of necklaces $N_1 \times \cdots \times N_k$ or equivalently with tuples $(N_1, \ldots, N_k)$ where $N_j \in \N_{\nu_j}$. Since
  \[ \Stab_{C_\nu}(N_1 \times \cdots \times N_k)
     = \prod_{j=1}^k \Stab_{C_{\nu_j}}(N_j), \]
we may group together $C_\nu$-orbits of $\W_n$ according to their
stabilizers as follows.

\begin{Definition}\label{def:NFD_nu_rho}
  Given $\nu = (\nu_1, \ldots, \nu_k)$ and
  $\rho = (\rho_1, \ldots, \rho_k)$, let
  \begin{align*}
    \NF_{\nu,\rho} &\coloneqq \NF_{\nu_1, \rho_1}
      \times \dots \times \NF_{\nu_k, \rho_k}, \\
    \NFD_{\nu,\rho} &\coloneqq \NFD_{\nu_1, \rho_1}
      \times \dots \times \NFD_{\nu_k, \rho_k}.
  \end{align*}
\end{Definition}
\noindent
The elements of $\NF_{\nu,\rho}$ all have the same
stabilizer, and the elements of $\NFD_{\nu,\rho}$
are precisely those whose stabilizer is contained in
the common stabilizer of elements of $\NF_{\nu,\rho}$.
We write $\rho \mid \nu$ to mean that $\rho_i \mid \nu_i$
for all $i=1, \ldots, r$. Note that $\NF_{\nu, \rho} \neq \varnothing$
if and only if $\rho \mid \nu$.

Given a group $G$ acting on a set $\W$ and a subgroup $H$ of $G$,
each $G$-orbit of $\W$ is partitioned into $H$-orbits. Consequently,
$C_\nu$-orbits of $\W_n$ are unions of $C$-orbits, which we exploit
as follows.

\begin{Lemma}\label{lem:C_Cnu.1}
  Let $\cO$ be a $C$-orbit of $\W_n$. Let
  $N_1 \times \cdots \times N_k$ be the $C_\nu$-orbit
  containing $\cO$ and suppose $N_1 \times \cdots \times N_k
  \in \NF_{\nu, \rho}$. Then
    \[ \#\cO = \lcm\left(\frac{\nu_1}{\rho_1}, \ldots,
       \frac{\nu_k}{\rho_k}\right), \]
  which depends only on $\nu$ and $\rho$. In particular,
    \[ \cO \in \OFD_{C, r} \qquad\text{ if and only if }\qquad
       \ell \mid r \cdot \lcm\left(\frac{\nu_1}{\rho_1},
       \ldots, \frac{\nu_k}{\rho_k}\right). \]

  \begin{proof}
    By assumption, $\freq(N_j) = \rho_j$ and $ N_j \in \N_{\nu_j} $, so
    $\#N_j = \nu_j/\rho_j$. It follows that $\cO$ is in
    bijection with the group generated by a permutation of
    cycle type $(\nu_1/\rho_1, \ldots, \nu_k/\rho_k)$, so
    that $\#\cO = \lcm(\nu_1/\rho_1, \ldots, \nu_k/\rho_k)$.
    The second claim follows by noting that
      \[ \cO \in \OFD_{C, r} \Leftrightarrow
         \freq(\cO) \mid r \Leftrightarrow
         (\ell/\#\cO) \mid r \Leftrightarrow
         \ell \mid r \cdot \#\cO. \]
  \end{proof}
\end{Lemma}

\begin{Lemma}\label{lem:OIRNFexpan}
  We have
    \[ \OFD_{C, r}^{\cont}(\mathbf{x})
          = \sum \f{\prod_{j = 1}^k \f{\nu_j}{\rho_j} }
              { \lcm \lp \f{\nu_1}{\rho_1}, \dots,
                \f{\nu_k}{\rho_k} \rp }
                \, \NF_{\nu, \rho}^{\cont}(\mathbf{x}), \]
  where the sum is over all $\rho$ such that $\rho \mid \nu$
  and $\ell \mid r \dd \lcm \lp \f{\nu_1}{\rho_1}, \dots,
  \f{\nu_k}{\rho_k} \rp$.

  \begin{proof}
    Consider the map
      \[
        \Omega : \OFD_{C, r} \to \bigsqcup \NF_{\nu,\rho}
      \]
    sending $ \cO \in \OFD_{C,r} $ to the $ C_\nu $-orbit containing
    $ \cO $,  where the union is over all $\rho$ such that
    $\rho \mid \nu $ and $\ell \mid r \dd \lcm \lp \f{\nu_1}{\rho_1},
    \dots, \f{\nu_k}{\rho_k} \rp $. By \Cref{lem:C_Cnu.1}, $ \Omega $
    does in fact map $  \OFD_{C, r} $ into this union, and $ \Omega $ is
    surjective. Also, each $ C_\nu $-orbit contained in
    $ \NF_{\nu, \rho} $ has size $ \lcm \lp \f{\nu_1}{\rho_1}, \dots,
    \f{\nu_k}{\rho_k} \rp $, and $ \# N_1 \times \dots \times N_k =
    \prod_{j=1}^k \frac{\nu_j}{\rho_j} $, so the fiber of each
    $ N_1 \times \cdots \times N_k \in \NF_{\nu, \rho} $ has size
     \[
        \# \Omega^{-1}(N_1 \x \dots \x N_k)
        = \frac{\prod_{j=1}^k \frac{\nu_j}{\rho_j}}
           {\lcm\left(\frac{\nu_1}{\rho_1}, \ldots,
           \frac{\nu_k}{\rho_k}\right)}.
      \]
    The result now follows from $ \Omega $ being content-preserving.
  \end{proof}
\end{Lemma}

In \Cref{sec:KW}, we used cyclic sieving to turn generating
functions involving $\NFD_{n, r}^{\cont}(\mathbf{x})$ into
Schur expansions. Thus our next goal is to turn the necklace
generating function in \Cref{lem:OIRNFexpan} into an
analogous generating function over
$\NFD_{\nu,\rho}^{\cont}(\mathbf{x})$. To accomplish this,
one could in principle use M\"obius inversion on the lattice of
stabilizers of $C_\nu$-orbits to convert from
$\NF_{\nu,\rho}^{\cont}(\mathbf{x})$ to
$\NFD_{\nu,\rho}^{\cont}(\mathbf{x})$. However, the following
argument is more direct.

\begin{Lemma}\label{lem:OIRNFDexpan}
  For $r=1, \ldots, n$,
    \[ \OFD_{C, r}^{\cont}(\mathbf{x})
       = \sum \NFD_{\nu, \tau}^{\cont}(\mathbf{x}), \]
  where the sum is over all $ k $-tuples of integers $\tau \in [\nu_1] \times \cdots \times [\nu_k]$
  such that $\sum_{j = 1}^k \f{\ell}{\nu_j} \tau_j \, \equiv_\ell \, r$.

  \begin{proof}
    We have
      \[ \NFD_{\nu, \tau}^{\cont}(\mathbf{x})
         = \sum_{\rho \mid \nu, \tau}
           \NF_{\nu, \rho}^{\cont}(\mathbf{x}), \]
where $\rho\mid\nu,\tau$ means $ \rho_j \mid \nu_j $ and $ \rho_j \mid \tau_j $ for all $ j $. Consequently,
    \begin{align*}
      \sum_{ \substack{ \tau \in [\nu_1] \times \dots \times [\nu_k] \\
             \sum_{j = 1}^k \f{\ell}{\nu_j} \tau_j \, \equiv_\ell \, r } }
        \NFD_{\nu, \tau}^{\cont}(\mathbf{x})
      = \sum_{ \rho \mid \nu } c_{\nu, \rho}^r
        \NF_{\nu, \rho}^{\cont}(\mathbf{x})
    \end{align*}
    where
    \begin{align*}
      c_{\nu,\rho}^r
        \coloneqq \#\left\{ \tau \in [\nu_1] \times \cdots \times [\nu_k]
             : \rho \mid \tau,
               \sum_{j=1}^k \frac{\ell}{\nu_j} \tau_j \equiv_\ell r
           \right\}.
    \end{align*}
    Since $\rho_j \mid \nu_j$ and $\rho_j \mid \tau_j$, write
    $ \gamma_j \coloneqq \frac{\nu_j}{\rho_j} \in \bZ_{\ge 1}$ and
    $ \de_j \coloneqq \frac{\tau_j}{\rho_j} \in \bZ_{\ge 1} $. Then,
      \[ \sum_{j=1}^k \frac{\ell}{\nu_j} \tau_j
         = \sum_{j=1}^k \frac{\ell}{\gamma_j} \de_j, \]
    so
      \[ c_{\nu, \rho}^r = \#\left\{ \de \in [\gamma_1] \times \cdots
           \times [\gamma_k] :
           \sum_{j=1}^k \frac{\ell}{\gamma_j} \de_j \equiv_\ell r \right\}. \]
    Defining a group homomorphism
    \begin{align*}
      \phi \colon \prod_{i=1}^k \bZ/\gamma_j &\to \bZ/\ell \\
      (\de_1, \ldots, \de_k) &\mapsto \sum_{j=1}^k \frac{\ell}{\gamma_j}\de_j,
    \end{align*}
    we now have $c_{\nu,\rho}^r = \#\phi^{-1}(r)$. Since
    $\frac{\ell}{\gamma_1}\bZ + \cdots + \frac{\ell}{\gamma_k}\bZ
    = \frac{\ell}{\lcm\lp\gamma_1, \ldots, \gamma_k\rp}\bZ$, it follows that
    \begin{align*}
      \im\phi &= \{r \in \bZ/\ell : \ell \mid r
        \cdot \lcm\lp\gamma_1, \ldots, \gamma_k\rp\}, \\
      \#\im\phi &= \lcm(\gamma_1, \ldots, \gamma_k).
    \end{align*}
    For $r \in \im \phi$, we then have
      \[ c_{\nu,\rho}^r = \#\phi^{-1}(r) = \#\ker \phi
         = \frac{\gamma_1 \cdots \gamma_k}
                {\lcm(\gamma_1, \ldots, \gamma_k)}. \]
    The result follows from \Cref{lem:OIRNFexpan}.
  \end{proof}
\end{Lemma}

Our next goal is to convert the necklace expansion in
\Cref{lem:OIRNFDexpan} into a Schur expansion. Recalling
from \Cref{sec:KW} that $ \M_{n, r} \coloneqq \{ w \in \W_n : \maj_n(w) = r \} $, \Cref{lem:NFD_F} tells us
\begin{equation}\label{eq:NFD_nutau}
  \NFD_{\nu,\tau}^{\cont}(\mathbf{x})
    = \prod_{j=1}^k \NFD_{\nu_j, \tau_j}^{\cont}(\mathbf{x})
    = \prod_{j=1}^k \M_{\nu_j, \tau_j}^{\cont}(\mathbf{x}).
\end{equation}
Interpreting the right-hand side of \eqref{eq:NFD_nutau} in terms of
words and comparing with the indexing set in
\Cref{lem:OIRNFDexpan} motivates the following
variations on the major index.

\begin{Definition}\label{def:majtuplestats}
  Suppose $\nu \vDash n$, $\tau \in [\nu_1] \times \cdots \times
  [\nu_k]$, and $\ell = \lcm(\nu_1, \ldots, \nu_k)$.
  Let $\bfmaj_\nu \colon \W_n \to [\nu_1] \times \cdots \times [\nu_k]$
  be defined as follows. For $w \in \W_n$, write
  $w = w^1 \cdots w^k$ where each $w^j$ is a word in $\W_{\nu_j}$. Set
    \[ \bfmaj_\nu(w) \coloneqq
       (\maj_{\nu_1}(w^1), \ldots, \maj_{\nu_k}(w^k)). \]
  Furthermore, let $\maj_\nu \colon \W_n \to [\ell]$ be defined by
    \[ \maj_\nu(w) \coloneqq
       \sum_{j=1}^k \frac{\ell}{\nu_j} \, \bfmaj_\nu(w)_j
       \qquad (\text{mod }\ell). \]
  Consequently, we have $\maj_{(n)} = \maj_n$.
  Note that both $\bfmaj_\nu$ and $\maj_\nu$ are functions of
  $\Des(w)$. We may thus define both $\bfmaj_\nu$ and
  $\maj_\nu$ on
  $Q \in \SYT(n)$ using only $\Des(Q)$ in the same way. Equivalently,
  we may set $\bfmaj_\nu(Q) \coloneqq \bfmaj_\nu(w)$ and
  $\maj_\nu(Q) \coloneqq \maj_\nu(w)$ for any $w$ such that $Q = Q(w)$.
\end{Definition}

\begin{Example}
  Let $\nu=(5, 3, 3)$ and $w=44121361631$, so that
  $\ell=15$, $w_1=44121$, $w_2=361$, and $w_3=631$.
  We have
    \[ \bfmaj_\nu(w) = (\maj_5(w_1), \maj_3(w_2), \maj_3(w_3))
       = (1, 2, 3) \]
  and, hence, $\maj_\nu(w) = \frac{15}{5} \cdot 1 + \frac{15}{3}
  \cdot 2 + \frac{15}{3} \cdot 3 = 13 \text{ $($mod $15)$}$.
\end{Example}

\begin{Definition}
  Suppose $\nu \vDash n$, $\tau \in [\nu_1] \times \cdots \times
  [\nu_k]$. Let
  \begin{align*}
    \M_{\nu,\tau} &\coloneqq \{w \in \W_n : \bfmaj_\nu(w) = \tau\}, \\
  \end{align*}
\end{Definition}

\begin{Theorem}{\cite[Theorem~3.3]{MR1023791}}\label{thm:Stembridge}
  Let $C$ be a cyclic subgroup of $S_n$ generated by an element of cycle type $\nu = (\nu_1, \dots, \nu_k) $, and let $ \ell = \lcm(\nu_1, \dots, \nu_k) $. We have
    \[ \sum_{r=1}^\ell \Ch\left(\chi^r\ind_C^{S_n}\right) q^r
       = \W_n^{\cont,\maj_\nu}(\mathbf{x}; q)
       = \sum_{\substack{\lambda \vdash n\\r \in [\ell]}}
         a_{\lambda, r}^\nu
         s_\lambda(\mathbf{x}) q^r \]
  where $a_{\lambda, r}^\nu \coloneqq \#\{Q \in \SYT(\lambda)
  : \maj_\nu(Q) = r\}$. In particular, the multiplicity of
  $S^\lambda$ in $\chi^r\ind_C^{S_n}$ is $a_{\lambda,r}^\nu$.

  \begin{proof}
    From the definition of $\bfmaj_\nu$ and \Cref{def:NFD_nu_rho}, we have
    \begin{equation}\label{eq:NFD_W.1}
      \NFD_{\nu,\tau}^{\cont}(\mathbf{x})
        = \M_{\nu,\tau}^{\cont}(\mathbf{x}).
    \end{equation}
    Using \Cref{thm:basis_cyclic} and \Cref{lem:OIRNFDexpan}, we
    then have
    \begin{align*}
      \sum_{r=1}^\ell \Ch\lp\chi^r\ind_C^{S_n}\rp q^r
        &= \sum_{r=1}^\ell \; \sum_{ \substack{ \tau \in [\nu_1] \times
              \dots \times [\nu_k] \\
              \sum_{j = 1}^k \f{\ell}{\nu_j} \tau_j \, \equiv_\ell \, r } }
            \NFD_{\nu, \tau}^{\cont}(\mathbf{x}) \; q^r \\
        &= \sum_{r=1}^\ell \; \sum_{ \substack{ \tau \in [\nu_1] \times
              \dots \times [\nu_k] \\
              \sum_{j = 1}^k \f{\ell}{\nu_j} \tau_j \, \equiv_\ell \, r } }
            \M_{\nu,\tau}^{\cont}(\mathbf{x}) \; q^r \\
        &= \sum_{r=1}^\ell
          \{ w \in \W_n : \maj_\nu(w) = r \}^{\cont}(\bfx) \; q^r \\
        &= \W_n^{\cont, \maj_\nu}(\mathbf{x}; q).
    \end{align*}
    Since $ \maj_\nu(w)$ depends only on $\Des(w)$, we can apply the RSK bijection again through \Cref{lem:RSK_Des} to get
      \[ \W_n^{\cont, \maj_\nu}(\mathbf{x}; q)
         = \sum_{\substack{\lambda \vdash n\\r \in [\ell]}}
           a_{\lambda, r}^\nu
           s_\lambda(\mathbf{x}) q^r. \]
  \end{proof}
\end{Theorem}

\begin{Remark}
  Stembridge showed the equality of the first and third terms
  in \Cref{thm:Stembridge} using the skew analogue of
  \eqref{eq:SYT_evals} and branching rules along
  Young subgroups of $S_n$.
  By contrast, $\W_n^{\cont, \maj_\nu}(\mathbf{x}; q)$
  played a key role in our approach.
\end{Remark}

Since the isomorphism type of $\chi^r\ind_C^{S_n}$, or equivalently
the Schur expansion of $\OFD_{C, r}^{\cont}(\mathbf{x})$, depends only on
$\nu $, the cycle type of a generator of $ C $, and $ \gcd(\ell, r) $, we have the following 
generalization of \Cref{cor:alamrsym}.

\begin{Corollary}\label{cor:alamrsym_nu}
  For all $ n \ge 1 $ and $ \lam, \nu \vdash n $, we have
  $ a_{\lam,r}^\nu = a_{\lam, \gcd(\ell,r)}^\nu $,
  where $\ell = \lcm(\nu_1, \nu_2, \ldots)$.
\end{Corollary}

For use in the next section, we record the Schur expansion of
$ \M_{\nu, \tau}^{\cont}(\mathbf{x}) $. The proof is analogous to the
last step of the proof of \Cref{thm:Stembridge} using \Cref{lem:RSK_Des}.

\begin{Corollary}\label{cor:M_schur}
  If $\nu,\tau \vDash n$, then
    \[ \M_{\nu, \tau}^{\cont}(\mathbf{x})
       = \sum_{\lambda \vdash n} \mathbf{a}_{\lambda,\tau}^\nu
         s_\lambda(\mathbf{x})  \]
  where
    \[ \mathbf{a}_{\lambda,\tau}^\nu \coloneqq
       \#\{Q \in \SYT(\lambda) : \bfmaj_\nu(Q) = \tau\}. \]
\end{Corollary}

We also have a corresponding symmetry result. Contrast it with
\Cref{cor:alamrsym}.

\begin{Corollary}\label{cor:permsym}
  Suppose $ \nu = (\nu_1, \dots, \nu_k) $ is the cycle type of
  some $\sigma \in S_n$, $ \tau \in [\nu_1] \times \cdots \times
  [\nu_k]$, $ \pi \in S_k $, and $ \lam \vdash n $. Then,
  $ \mathbf{a}_{\lam,\tau}^{\nu} =
  \mathbf{a}_{\lam, \pi \dd \tau}^{\pi \dd \nu} $.

  \begin{proof}
    Since reordering does not affect contents, we have
      \[ \NFD_{\nu, \tau}^{\cont}(\mathbf{x})
         = \NFD_{\pi \cdot \nu, \pi \cdot \tau}^{\cont}
           (\mathbf{x}). \]
    Now apply \Cref{cor:M_schur} and equate coefficients of
    $s_\lambda(\mathbf{x})$.
  \end{proof}
\end{Corollary}

\section{Inducing 1-dimensional Representations from \texorpdfstring{$ C_a \wr S_b $}{Ca wreath Sb} to \texorpdfstring{$ S_{ab} $}{Sab}}\label{sec:HLM}
We next apply the approach of \Cref{sec:KW} and \Cref{sec:Cn_branching}
to prove a generalization of a formula due to Schocker
\cite{MR1984625} for the Schur expansion of $\cL_{(a^b)}$.
In particular, we give Schur expansions of the characteristics of
  \[ \cL_{(a^b)}^{r, 1}
       \coloneqq \chi^{r, 1}\ind_{C_a \wr S_b}^{S_{ab}}
     \qquad \text{ and } \qquad
     \cL_{(a^b)}^{r, \epsilon}
       \coloneqq \chi^{r, \epsilon}\ind_{C_a \wr S_b}^{S_{ab}}. \]
Note that $\Ch \cL_{(a^b)} = \Ch \cL_{(a^b)}^{1, 1}$ by
\Cref{cor:higher_schur_weyl}.

The argument in \Cref{cor:higher_schur_weyl} and the fact that $ \Ch(\epsilon_b) = e_b(\bfx) $ immediately
yield the following more general result, which also
follows from an appropriate modification of \Cref{thm:basis_onerow}.

\begin{Lemma}\label{lem:setmultisetGFs}
  We have
    \[ \cL_{(a^b)}^{r, 1}
       = \mch{\NFD_{a, r}}{b}^{\cont}(\mathbf{x})
       \qquad \text{ and } \qquad
       \cL_{(a^b)}^{r, \epsilon}
       = \binom{\NFD_{a, r}}{b}^{\cont}(\mathbf{x}). \]
\end{Lemma}

Our first goal is to manipulate the necklace generating functions in
\Cref{lem:setmultisetGFs} in such a way that we may apply cyclic
sieving. We use Burnside's lemma and a sign-reversing involution
to unravel these multiset and subset generating functions, respectively.

\begin{Lemma}\label{lem:schocker_burnside}
  We have
    \[ \mch{\NFD_{a, r}}{b}^{\cont}(x_1, x_2, \ldots)
       = \sum_{\nu \vdash b} \frac{1}{z_\nu} \prod_{j=1}^{\ell(\nu)}
         \NFD_{a, r}^{\cont}(x_1^{\nu_j}, x_2^{\nu_j}, \ldots). \]

  \begin{proof}
    Multisets of $b$ necklaces from $\NFD_{a, r}$ can be thought of as
    $S_b$-orbits of length-$b$ tuples $(N_1, \ldots, N_b)$ of necklaces
    $N_i \in \NFD_{a, r}$ under the natural $S_b$-action. The tuples
    $(N_1, \ldots, N_b)$ fixed by an element $\si \in S_b$ are those
    tuples which are constant on blocks corresponding to cycles of $\si$.
    It follows that if $\si$ has cycle type $\nu \vdash b$,
    \begin{equation}\label{eq:contgf_fixednecklaces}
      \{ T \in \NFD_{a,r}^b : \si \dd T = T \}^{\cont}
      (x_1, x_2, \ldots) = \prod_{j=1}^{\ell(\nu)}
      \NFD_{a, r}^{\cont}(x_1^{\nu_j}, x_2^{\nu_j}, \ldots).
    \end{equation}
    By Burnside's lemma, we may count $S_b$-orbits of necklaces
    $(N_1, \ldots, N_b)$ of fixed content by averaging the number
    of $\si$-fixed tuples of fixed content over all $\si \in S_b$. The
    result follows by grouping together permutations of a given
    cycle type.
  \end{proof}
\end{Lemma}

\begin{Lemma}\label{lem:schocker_burnside.eps}
  We have
    \[ \binom{\NFD_{a, r}}{b}^{\cont}(\bfx)
       = \sum_{\nu \vdash b} \frac{(-1)^{b - \ell(\nu)}}{z_\nu}
         \prod_{j=1}^{\ell(\nu)}
         \NFD_{a, r}^{\cont}(x_1^{\nu_j}, x_2^{\nu_j}, \ldots). \]
\end{Lemma}

\begin{proof}
  Multiplying both sides by $b!$, using
  \eqref{eq:contgf_fixednecklaces} and the fact
  $ \sgn(\si) = (-1)^{b - \ell(\nu)} $ for $ \si \in S_b $
  with cycle type $ \nu $, the result is equivalent to
  \begin{equation}\label{eq:fixednecklacesSRI}
    \begin{aligned}
      \{ (N_1, \dots, N_b) & \in \NFD_{a, r}^b
        : (N_1, \dots, N_b) \tx{ are distinct} \}^{\cont}(\bfx) \\
      &= \sum_{\si \in S_b} \sgn(\si)
         \{ T \in \NFD_{a,r}^b : \si \dd T = T \}^{\cont}(\bfx).
    \end{aligned}
  \end{equation}
  On the right-hand side of \eqref{eq:fixednecklacesSRI}, each $b$-tuple
  $(N_1, \ldots, N_b)$ is counted
  \[
      \wgt(N_1, \dots, N_b) \coloneqq
      \sum_{\substack{\si \in S_b \\
        \tx{s.t. } \si \dd (N_1, \dots, N_b) = (N_1, \dots, N_b)} }
          \sgn(\si)
  \]
  times. If $N_1, \ldots, N_b$ are distinct, then only $ \si =\id$
  contributes, so $ \wgt(N_1, \dots, N_b) = 1 $. If
  $ N_1, \ldots, N_b$ are not distinct, then without loss of
  generality, suppose $N_1 = N_2$. Then, modifying the cycle(s) 
  containing $1$ and $2$ as in
      \[ (1\ \cdots)(2\ \cdots) \leftrightarrow (1\ \cdots\ 2\ \cdots) \]
  gives a sign-reversing involution on
  $ \{ \si \in S_b : \si \dd (N_1, \dots, N_b) = (N_1, \dots, N_b) \} $,
  meaning $ \wgt(N_1, \dots, N_b) = 0 $. This proves
  \eqref{eq:fixednecklacesSRI}.
\end{proof}

\begin{Remark}
  Using standard properties of plethysm (see
  e.g.~\cite[\S I.8]{MR1354144})) and the power-sum expansions
  of $e_b$ and $h_b$ (see \cite[(7.22)-(7.23)]{MR1676282}),
  \Cref{lem:schocker_burnside} and \Cref{lem:schocker_burnside.eps}
  are equivalent to
  \begin{align}\label{eq:schocker_p.1}
    \Ch\cL_{(a^b)}^{r, 1}
     &= h_b[\Ch \chi^r\ind_{C_a}^{S_a}]
      = \sum_{\nu \vdash b} \frac{1}{z_\nu}
       p_\nu[\Ch \chi^r\ind_{C_a}^{S_a}], \\
    \label{eq:schocker_p.2}
    \Ch\cL_{(a^b)}^{r, \epsilon}
      &= e_b[\Ch \chi^r\ind_{C_a}^{S_a}]
       = \sum_{\nu \vdash b}
       \frac{(-1)^{b-\ell(\nu)}}{z_\nu}
       p_\nu[\Ch \chi^r\ind_{C_a}^{S_a}].
  \end{align}
  Consequently, one may replace the combinatorial manipulations
  in \Cref{lem:schocker_burnside} and \Cref{lem:schocker_burnside.eps}
  with symmetric function manipulations. In the next section,
  we will prove \Cref{thm:casbchars}, which generalizes the first
  equalities in \eqref{eq:schocker_p.1} and \eqref{eq:schocker_p.2}.
\end{Remark}

\begin{Remark}
  Let $\omega$ be the involution on the algebra of symmetric
  functions defined by $\omega(s_\lam(\bfx)) = s_{\lam'}(\bfx)$
  where $\lambda'$ is the \textit{conjugate} of $\lambda$,
  obtained by reflecting $\lambda$ through the line $y=-x$.
  One may show in a variety of ways that
  \begin{align}\label{eq:omegaind}
    \omega\lp\Ch \chi^r\ind_{C_n}^{S_n}\rp
      = \Ch \chi^s \ind_{C_n}^{S_n}
        \qquad\text{where}\qquad s = \binom{n}{2} - r.
  \end{align}

  For instance, we can prove \eqref{eq:omegaind} using \Cref{thm:KW}
  as follows. Since conjugation $ Q \mapsto Q' $ satisfies
  $ \Des(Q') = [n - 1] \sm \Des(Q) $, we have
  \begin{equation}
    \begin{aligned}
      a_{\lam',r}
        &= \# \{ Q \in \SYT(\lam') : \maj(Q) \equiv_n r \} \\
        &= \# \left\{ Q' \in \SYT(\lam) : \maj(Q')
           \equiv_n \ch{n}{2} - r \right\} = a_{\lam, \ch{n}{2} - r}.
    \end{aligned}
  \end{equation}
  Therefore, by \Cref{thm:KW}, letting $ s = \ch{n}{2} - r $,
  \begin{align*}
    \omega\lp\Ch \chi^r\ind_{C_n}^{S_n} \rp
      = \omega \lp \sum_{\lam \vdash n} a_{\lam, r} s_{\lam}(\bfx) \rp
      = \sum_{\lam \vdash n} a_{\lam', r} s_{\lam}(\bfx)
      = \sum_{\lam \vdash n} a_{\lam, s} s_{\lam}(\bfx)
      = \chi^s \ind_{C_n}^{S_n}.
  \end{align*}

  From the symmetry result \Cref{cor:alamrsym},
  it follows that $\Ch \chi^r\ind_{C_n}^{S_n}$ is
  fixed under $\omega$ when $n$ is odd. When $n$ is even,
  $\Ch \chi^r\ind_{C_n}^{S_n}$ may or may not be fixed.
  For instance, when $r=1$, we find
    \[
      \omega\lp\Ch\cL_{n}\rp =
      \omega\lp\Ch\chi^1\ind_{S_n}^{C_n}\rp =
      \begin{cases}
        \Ch\cL_{n}^{(2)}
          = \Ch\chi^2\ind_{S_n}^{C_n}
          & \text{if $n/2$ is odd,} \\
        \Ch\cL_{n}
          = \Ch\chi^1\ind_{S_n}^{C_n}
          & \text{otherwise.}
      \end{cases}
    \]
  Here $\cL_n^{(2)}$ is the deformation of
  $\cL_n$ recently studied by Sundaram \cite{arX1803.09368}.
  Further standard properties of plethysm together
  with \eqref{eq:schocker_p.1} and \eqref{eq:schocker_p.2} give
  \begin{align*}
    &\ \omega\lp\Ch\cL_{(a^b)}^{r, 1}\rp
    = \Ch\cL_{(a^b)}^{r, \epsilon}
      \qquad\text{for $a$ odd, and} \\
    &\left.
    \begin{aligned}
    \omega\lp\Ch\cL_{(a^b)}^{r, 1}\rp
      &= \Ch\cL_{(a^b)}^{s, 1} \\
    \omega\lp\Ch\cL_{(a^b)}^{r, \epsilon}\rp
      &= \Ch\cL_{(a^b)}^{s, \epsilon}
    \end{aligned}
    \qquad\right\}\quad \text{for $a$ even, where $s=\binom{a}{2} - r$.}
  \end{align*}
  Consequently, one may obtain the Schur expansion of
  $\Ch \cL_{(a^b)}^{r, \epsilon}$ from
  the Schur expansion of $\Ch \cL_{(a^b)}^{r, 1}$ simply by
  applying the $\omega$ map if and only if $a$ is odd.
  When $a$ is even, these two cases are more fundamentally different.
\end{Remark}

Next, we convert
$ \NFD_{a, r}^{\cont}(x_1^{\nu_j}, x_2^{\nu_j}, \ldots) $ into a
linear combination of $\NF_{k, s}^{\cont}(\mathbf{x})$'s and then
apply Mobius inversion to convert to a linear combination of
$\NFD_{k, s}^{\cont}(\mathbf{x})$'s. We will need the following
variation on the number-theoretic M\"obius function $\mu$.

\begin{Definition}\label{def:mobiusf}
  Suppose $d \mid e$ and $f \mid e$. Set
    \[ \mu_f(d, e)
       \coloneqq \sum_{\substack{g \\
         \text{s.t. } \lcm(f, d) \mid g \mid e}}
            \mu\lp\frac{g}{f}\rp. \]
\end{Definition}

This expression simplifies considerably as follows. Let
$ \rad(m) $ denote the squarefree positive integer with
the same prime divisors as $m$.

\begin{Lemma}\label{lem:mobiusf}
  Suppose $d \mid e$ and $f \mid e$. Then
  \begin{align*}
    \mu_f(d, e)
      = \begin{cases}
          \mu\lp\frac{\lcm(f, d)}{f}\rp
            & \text{if }\rad\lp\frac{e}{f}\rp
                = \rad\lp\frac{\lcm(f, d)}{f}\rp
                = \frac{\lcm(f, d)}{f}, \\
          0 & \text{otherwise.}
        \end{cases}
  \end{align*}

  \begin{proof}
    We see
    \begin{align*}
      \mu_f(d,e) = \sum_{\substack{g \\ \text{s.t. } \lcm(f, d) \mid g \mid e}}
          \mu\lp\frac{g}{f}\rp
        = \sum_{\substack{h \\ \text{s.t. }
                \frac{\lcm(f, d)}{f} \mid h \mid \frac{e}{f}}} 
	\mu(h)
    \end{align*}
Since $ \mu(h) \ne 0 $ only when $ h $ is radical, $ \mu(h) \ne 0 $ only when $ \frac{\lcm(f, d)}{f} $ is radical and $ h \mid \rad(e/f) $. Restricting to this case, we can write $ \rad(e/f) = k \lcm(f, d)/f $ for some integer $ k $. Since $ k $ and $\lcm(f, d)/f$ must be relatively prime, we have
    \begin{align*}
      \mu_f(d, e)
        &= \sum_{\substack{h \\ \text{s.t. }
                \frac{\lcm(f, d)}{f} \mid h \mid
                \frac{k \cdot \lcm(f, d)}{f}}}
           \mu(h) \\
        &= \sum_{s \mid k}
           \mu\lp\lp\frac{\lcm(f, d)}{f}\rp s \rp \\
        &= \mu\lp\frac{\lcm(f, d)}{f}\rp \sum_{s \mid k} \mu(s) \\
        &= \begin{cases}
          \mu\lp\frac{\lcm(f, d)}{f}\rp
            & \text{if } k = 1, \\
          0 & \text{otherwise},
        \end{cases}
    \end{align*}
    giving the result.
  \end{proof}
\end{Lemma}

\begin{Lemma}\label{lem:NFD_pleth}
  We have
    \[ \NFD_{a, r}^{\cont}(x_1^k, x_2^k, \ldots)
       = \sum_{s \mid rk} \mu_s(k, rk) \NFD_{ak, s}^{\cont}(x_1, x_2, \ldots). \]

  \begin{proof}
    The left-hand side is the content generating function for
    $k$-tuples of length $a$ necklaces with frequency dividing
    $ r $ of the form $ (N, \dots, N) $, repeating the same necklace
    $ k $ times. By concatenation, we may equivalently view such
    tuples as length $ak$ necklaces whose frequency $f$ satisfies
    $k \mid f \mid rk$. Consequently,
    \begin{align}\label{eq:freqsum}
      \NFD_{a, r}^{\cont}(x_1^k, x_2^k, \ldots)
        = \sum_{\substack{f \\\text{s.t. }k \mid f \mid rk}}
          \NF_{ak, f}^{\cont}(x_1, x_2, \ldots),
    \end{align}
    recalling $ \NF_{n,f} \coloneqq \{ N \in \N_n : \freq(N) = f \} $.
    M\"obius inversion on the identity
    $ \NFD_{ak, f}^{\cont}(\mathbf{x})
    = \sum_{s \mid f} \NF_{ak, s}^{\cont}(\mathbf{x}) $ gives
    \begin{align}\label{eq:Mobiusinv}
      \NF_{ak, f}^{\cont}(\mathbf{x})
      = \sum_{s \mid f} \mu \lp \frac{f}{s} \rp
        \NFD_{ak, s}^{\cont}(\bfx).
    \end{align}
    Thus, by \eqref{eq:Mobiusinv}, \eqref{eq:freqsum} becomes
    \begin{align*}
      \NFD_{a, r}^{\cont}(x_1^k, x_2^k, \ldots)
        &= \sum_{\substack{f \\ \text{s.t. }k \mid f \mid rk}}
           \sum_{s \mid f} \mu\lp\frac{f}{s}\rp
             \NFD_{ak, s}^{\cont}(x_1, x_2, \ldots) \\
        &= \sum_{s \mid rk} \lp
           \sum_{ \substack{f \\ \text{s.t. } \lcm(k,s) \mid f \mid rk}}
             \mu\lp \frac{f}{s}\rp \rp
             \NFD_{ak, s}^{\cont}(x_1, x_2, \ldots) \\
        &= \sum_{s \mid rk} \mu_s(k, rk)
             \NFD_{ak, s}^{\cont}(x_1, x_2, \ldots).
    \end{align*}
    by \Cref{def:mobiusf}.
  \end{proof}
\end{Lemma}

\begin{Notation}
  Given a sequence $\nu = (\nu_1, \ldots, \nu_k) \in \bZ_{\geq 1}^k$
  and an integer $r \in \bZ_{\geq 1}$, let
    \[ r \ast \nu \coloneqq (r\nu_1, \ldots, r\nu_k). \]
  Given another sequence $\tau = (\tau_1, \ldots, \tau_k)$,
  recall that $\tau \mid \nu$ means $\tau_j \mid \nu_j$ for all $j$.
  Further recall
  \[
      \NFD_{\nu, \tau}
      = \NFD_{\nu_j, \tau_j} \times \dots \times \NFD_{\nu_k, \tau_k}
  \]
  from \Cref{def:NFD_nu_rho}. Finally, extend $\mu_f(d, e)$ to
  sequences multiplicatively:
    \[ \mu_{(f_1, \ldots, f_k)}((d_1, \ldots, d_k), (e_1, \ldots, e_k))
       \coloneqq \prod_{j=1}^k \mu_{f_j}(d_j, e_j). \]
\end{Notation}

\begin{Corollary}\label{cor:NFD_mult_NFD}
  We have
  \begin{align*}
    \Ch \cL_{(a^b)}^{r, 1}
      &= \sum_{\nu \vdash b} \frac{1}{z_\nu} \sum_{\tau \mid r \ast \nu}
         \mu_\tau(\nu, r\ast\nu)
         \NFD_{a\ast\nu, \tau}^{\cont}(\mathbf{x}), \\
    \Ch \cL_{(a^b)}^{r, \epsilon}
      &= \sum_{\nu \vdash b} \frac{(-1)^{b - \ell(\nu)}}{z_\nu}
         \sum_{\tau \mid r \ast \nu}
         \mu_\tau(\nu, r\ast\nu)
         \NFD_{a\ast\nu, \tau}^{\cont}(\mathbf{x}).
  \end{align*}

  \begin{proof}
    Combine \Cref{lem:setmultisetGFs}, \Cref{lem:schocker_burnside} or
    \Cref{lem:schocker_burnside.eps}, and \Cref{lem:NFD_pleth}.
  \end{proof}
\end{Corollary}

We may now state and generalize Schocker's formula for
$\Ch \cL_{(a^b)} = \Ch \cL_{(a^b)}^{1, 1}$.

\begin{Theorem}[See {\cite[Thm.~3.1]{MR1984625}}]\label{thm:GeneralizedSchocker}
  For all $a, b \geq 1$ and $r=1, \ldots, a$, we have
  \begin{align*}
    \Ch \cL_{(a^b)}^{r, 1}
      &= \sum_{\lambda \vdash ab} \lp\sum_{\nu \vdash b} \frac{1}{z_\nu}
         \sum_{\tau \mid r \ast \nu} \mu_\tau(\nu, r \ast \nu)
           \mathbf{a}_{\lambda, \tau}^{a \ast \nu}\rp
           s_\lambda(\mathbf{x}) \qquad \text{and} \qquad \\
    \Ch \cL_{(a^b)}^{r, \epsilon}
      &= \sum_{\lambda \vdash ab} \lp\sum_{\nu \vdash b}
         \frac{(-1)^{b - \ell(\nu)}}{z_\nu}
         \sum_{\tau \mid r \ast \nu} \mu_\tau(\nu, r \ast \nu)
           \mathbf{a}_{\lambda, \tau}^{a \ast \nu}\rp
           s_\lambda(\mathbf{x}),
  \end{align*}
  where, recalling the definition of $ \mathbf{maj}_{a \ast \nu} $ from
  \Cref{def:majtuplestats},
    \[ \mathbf{a}_{\lambda, \tau}^{a \ast \nu} \coloneqq
       \#\{Q \in \SYT(\lambda) : \mathbf{maj}_{a \ast \nu}(Q)
       = \tau\}. \]

  \begin{proof}
    Combine \Cref{cor:M_schur} and \Cref{cor:NFD_mult_NFD}.
  \end{proof}
\end{Theorem}

\begin{Remark}
  Schocker's approach to \cite[Thm.~3.1]{MR1984625} uses J\"ollenbeck's
  non-commutative character theory and involved manipulations with
  Klyachko's idempotents and Ramanujan sums. Much of Schocker's argument
  generalizes immediately to all $r$. The argument presented above
  is comparatively self-contained and direct. Two perhaps mysterious
  aspects of the formula, the appearance of M\"obius functions and
  the average over $S_b$, arose naturally from Burnside's lemma and
  a change of basis using M\"obius inversion. Our argument uses explicit
  bijections at each step except for the appeal to Burnside's lemma and
  the use of \Cref{lem:schocker_burnside.eps}.
\end{Remark}

\section{Higher Lie Modules and Branching Rules}\label{sec:mash}
The argument in \Cref{sec:KW} solves Thrall's problem for
$\lambda=(n)$ by considering all branching rules for
$C_n \hookrightarrow S_n$ simultaneously and using cyclic sieving
and RSK to convert from the monomial to the Schur basis.
We now turn to analogous considerations for the higher Lie modules
and more generally branching rules for
$C_a \wr S_b \hookrightarrow S_{ab}$. We give an analogue of
the $\flex$ statistic and the monomial basis expansion for
such branching rules from \Cref{ssec:background_KW}. We then show
how to convert from the monomial to the Schur basis assuming
the existence of a certain statistic on words we call $\mash$
which interpolates between $\maj_n$ and the shape under RSK.

We now recall and prove \Cref{thm:casbchars} from the introduction, after introducing some notation.

\begin{Definition} Fix integers $ a, b \ge 1$. Define
\[
  \PP_{a}^{b} \coloneqq \left\{ \ul = (\lam^{(1)}, \dots, \lam^{(a)}) : \lam^{(1)}, \dots, \lam^{(a)} \tx{ are partitions }, \sum_{r = 1}^a |\lam^{(r)}| = b \right\},
\]
  which indexes the irreducible $C_a \wr S_b$-representations by
  \Cref{thm:CaSb_irreps}.
\end{Definition}

\begin{Theorem*}
  For all $ a, b \ge 1$ and $\ul = (\lam^{(1)}, \dots, \lam^{(a)}) \in \PP_a^b $, we have
  \[ \Ch S^{\ul}\ind_{C_a \wr S_b}^{S_{ab}}
     = \prod_{r = 1}^{a}
     s_{\lambda^{(r)}}[\NFD_{a,r}^{\cont}(\mathbf{x})]. \]

  \begin{proof}[Proof of {\Cref{thm:casbchars}}]
  We have
  \begin{align*}
    S^{\ul}\ind_{C_a \wr S_b}^{S_{ab}}
      & = \left[\bigotimes_{r=1}^a (\chi^r_a \wr S^{\lambda^{(r)}})\right]
        \ind_{C_a \wr S_{\al(\ul)}}^{S_{ab}} \\
      &\cong \left[\bigotimes_{r=1}^a (\chi^r_a \wr S^{\lambda^{(r)}})\right]
        \ind_{C_a \wr S_{\al(\ul)}}^{S_{a \ast \al(\ul)}}
        \ind_{S_{a \ast \al(\ul)}}^{S_{ab}} \\
      &\cong \left[\bigotimes_{r=1}^a
        (\chi^r_a \wr S^{\lambda^{(r)}})
        \ind_{C_a \wr S_{|\lambda^{(r)}|}}^{S_{a|\lambda^{(r)}|}}\right]
        \ind_{S_{a \ast \al(\ul)}}^{S_{ab}} \\
      &\cong \left[\bigotimes_{r=1}^a
        (\chi^r_a \wr S^{\lambda^{(r)}})
        \ind_{C_a \wr S_{|\lambda^{(r)}|}}^{S_a \wr S_{|\lambda^{(r)}|}}
        \ind_{S_a \wr S_{|\lambda^{(r)}|}}^{S_{a|\lambda^{(r)}|}}\right]
        \ind_{S_{a \ast \al(\ul)}}^{S_{ab}} \\
      &\cong \left[\bigotimes_{r=1}^a
        (\chi^r_a\ind_{C_a}^{S_a} \wr S^{\lambda^{(r)}})
        \ind_{S_a \wr S_{|\lambda^{(r)}|}}^{S_{a|\lambda^{(r)}|}}\right]
        \ind_{S_{a \ast \al(\ul)}}^{S_{ab}},
  \end{align*}
  where the first and third isomorphisms use transitivity of
  induction, the second isomorphism uses \Cref{lem:indtensor}, and
  the fourth isomorphism uses \Cref{lem:indwreath}.
  Consequently, using \eqref{eq:prodchar}, \eqref{eq:plethsymchar},
  and \Cref{thm:basis_onerow}, we have
  \begin{align*}
    \Ch S^{\ul}\ind_{C_a \wr S_b}^{S_{ab}}
      &= \prod_{r=1}^a \Ch \lp\chi^r_a\ind_{C_a}^{S_a} \wr S^{\lambda^{(r)}}\rp
         \ind_{S_a \wr S_{|\lambda^{(r)}|}}^{S_{a|\lambda^{(r)}|}} \\
      &= \prod_{r=1}^a (\Ch S^{\lambda^{(r)}})[\Ch \chi^r_a\ind_{C_a}^{S_a}] \\
      &= \prod_{r=1}^a s_{\lambda^{(r)}}[\NFD_{a, r}^{\cont}(\mathbf{x})].
  \end{align*}
  \end{proof}
\end{Theorem*}

Recall from \Cref{ssec:background_tableaux} that given a word $w$,
the shape of $w$, denoted $\sh(w)$, is the common shape of $P(w)$
and $Q(w)$ under RSK.

\begin{Definition}\label{def:flex_ab}
  Fix $a, b \ge 1 $. Construct statistics
  \[ \flex_a^b, \maj_a^b \colon \W_{ab}
     \to \PP_a^b \]
  as follows. Given $w \in \W_{ab}$, write $w = w^1 \cdots w^b $ where
  $ w^j \in \W_a $. In this way, consider $ w $ as a word of size
  $ b $ whose letters are in $ \W_a $. For each $r \in [a]$,
  let $w^{(r)}$ denote the subword of $w$ whose letters are those
  $w^j$ such that $\flex(w^j) = r$. Totally order $\W_a$
  lexicographically, so that RSK is well-defined
  for words with letters from $\W_a$. Set
    \[ \flex_a^b(w) \coloneqq (\sh(w^{(1)}), \ldots, \sh(w^{(a)})). \]
  Define $\maj_a^b$ in the same way but with $\flex$ replaced by
  $\maj_a$. Consequently, $\maj_n^1(w)$ is the $n$-tuple of
  partitions whose only non-empty entry is a single cell
  at position $\maj_n(w)$.
\end{Definition}

\begin{Example}
  Let $w = 212023101241$ and suppose $a=3$, $b=4$. Write
  $ w = (212)(023)(101)(241)$. The parenthesized terms have
  $ \flex $ statistics $ 2, 1, 2, 2 $ and
  $ \maj_3 $ statistics $ 1, 3, 1, 2 $, respectively.
  When computing $\flex_3^4(w)$, we then have
  $ w^{(1)} = (023), w^{(2)} = (212)(101)(241), w^{(3)} = \vn $. Since
  $ (101) <_{\lex} (212) <_{\lex} (241), $
  $\sh(w^{(2)}) = \sh(213) = (2, 1)$. Consequently,
    \[ \flex_3^4(212023101241) = ((1), (2, 1), \varnothing).
    \]
  When computing $\maj_3^4(w)$, we have
  $ w^{(1)} = (212)(101), w^{(2)} = (241), w^{(3)} = (023) $.
  Since $ (101) <_{\lex} (212) $, $ \sh(w^{(1)}) = \sh(21) = (1,1)$.
  Hence
    \[ \maj_3^4(212023101241) = ((1, 1), (1), (1)). \]
\end{Example}

We now recall and prove \Cref{thm:grfrob_flexab} from the
introduction.

\begin{Theorem*}
  Fix $ a, b \ge 1 $. We have
\begin{align*}
    \sum_{\ul \in \PP_a^b }
      \dim S^{\ul} \cdot
      \Ch \lp S^{\ul}\ind_{C_a \wr S_b}^{S_{ab}} \rp q^{\ul}
      &= \W_{ab}^{\cont,\flex_a^b}(\mathbf{x}; q) \\
      &= \W_{ab}^{\cont,\maj_a^b}(\mathbf{x}; q)
\end{align*}
  where the $ S^{\ul} $ are irreducible representations of
  $ C_a \wr S_b $ and the $q^{\ul}$ are independent indeterminates.
\end{Theorem*}

\begin{proof}[Proof of {\Cref{thm:grfrob_flexab}}.]
  Fix $\ul \in \PP_a^b $. For the left-hand side, using
  \Cref{thm:casbchars} and \eqref{eq:dimcasbirreps},
  \begin{align}\label{eq:charabGF}
    \dim S^{\ul} \cdot
      \Ch \lp S^{\ul}\ind_{C_a \wr S_b}^{S_{ab}} \rp
      = \ch{b}{\al(\ul)} \prod_{r = 1}^a \#\SYT(\lambda^{(r)}) \cdot
        s_{\lambda^{(r)}}[\NFD_{a,r}^{\cont}(\mathbf{x})].
  \end{align}
  For the right-hand side, we have
    \[ \left. \W_{ab}^{\cont, \flex_a^b}(\mathbf{x}; q)
       \right|_{q^{\ul}}
       = \{ w \in \W_{ab} :
            \flex_a^b(w) = \ul \}^{\cont}(\mathbf{x}). \]
  Say $ \al(\ul) = (\al_1, \dots, \al_{a}) $. In order for
  $ w \in \W_{ab} $ to have $ \flex_a^b(w) = \ul $,
  we must have $\sh(w^{(r)}) = \lambda^{(r)}$ for each
  $r \in [a]$. Recalling $\F_{a,r} \coloneqq
  \{w \in \W_a : \flex(w) = r\}$,
  we may thus choose each $w^{(r)} \in (\F_{a,r})^{\alpha_r}$
  with $\sh(w^{(r)}) = \lambda^{(r)}$ independently and then shuffle them
  in $\binom{b}{\alpha(\ul)}$ ways to form $w$. Consequently,
  \begin{equation}\label{eq:flexshuffle}
    \begin{aligned}
      \{ w \in \W_{ab} &: \flex_a^b(w) = \ul \}^{\cont}(\mathbf{x}) \\
        &= \binom{b}{\alpha(\ul)} \prod_{r=1}^a
           \{w^{(r)} \in (\F_{a,r})^{\alpha_r} : \sh(w^{(r)})
         = \lambda^{(r)}\}^{\cont}(\mathbf{x}).
    \end{aligned}
  \end{equation}
  The content generating function for words with a given shape
  $\mu \vdash n$ under RSK is given by
    \begin{align}
      \{ w \in \W_n : \sh(w) = \mu \}^{\cont}(\mathbf{x})
      = \# \SYT(\mu) \, s_\mu(\mathbf{x}),
    \end{align}
  since the number of possible $Q$-tableaux is $ \# \SYT(\mu) $ and
  the content generating function for $P$-tableaux is
  $ s_\lam(\mathbf{x}) $.
  Changing the alphabet from $\bZ_{\geq 1}$ to $ \F_{a,r} $
  and using \Cref{lem:NFD_F} gives
  \begin{equation}\label{eq:Flexalphabet}
    \begin{aligned}
      \{ w^{(r)} \in (\F_{a,r})^{\al_r} : \sh(w^r)
        = \lambda^{(r)} \}^{\cont}(\bfx) &
        = \# \SYT(\lambda^{(r)})
          s_{\lambda^{(r)}}[\F_{a,r}^{\cont}(\mathbf{x})] \\
       &= \# \SYT(\lambda^{(r)})
          s_{\lambda^{(r)}}[\NFD_{a,r}^{\cont}(\mathbf{x})].
    \end{aligned}
  \end{equation}
  The first equality in \Cref{thm:grfrob_flexab} now follows from
  combining \eqref{eq:flexshuffle} and \eqref{eq:Flexalphabet} with
  \eqref{eq:charabGF}. The second equality in \Cref{thm:grfrob_flexab}
  follows similarly.
\end{proof}

While \Cref{thm:grfrob_flexab} determines the monomial
expansion of the graded Frobenius series tracking branching
rules for $C_a \wr S_b \hookrightarrow S_{ab}$, we are
ultimately interested in the corresponding Schur expansion.
We next describe how the approach in the preceding sections
might be used to find this Schur expansion. The key properties
used in the proof of \Cref{thm:KW} converting from the
monomial basis to the Schur basis were that $\maj_n$ is equidistributed with
$\flex$ on each $\W_\alpha$ and $ \maj_n(w) $ depends
only on $ Q(w) $. In order to apply a similar argument for
$ \Ch( S^{\ul}\ind_{C_a \wr S_b}^{S_{ab}} ) $, we need a
statistic as follows.

\begin{Problem}\label{prob:mash}
  Fix $a, b  \ge 1 $. Find a statistic
    \[ \mash_a^b \colon \W_{ab} \to
       \PP_{a}^{b} \]
  with the following
  properties.
  \begin{enumerate}[(i)]
    \item For all $\alpha \vDash ab$, $\maj_a^b$ (or equivalently
      $\flex_a^b$) and $\mash_a^b$ are equidistributed on $ \W_\al $.
    \item If $ v, w \in \W_{ab} $ satisfy $ Q(v) = Q(w) $,
      then $\mash_a^b(v) = \mash_a^b(w) $.
  \end{enumerate}
\end{Problem}

Finding such a statistic $ \mash_a^b $ would determine
the Schur decomposition of $ \Ch( S^{\ul}
\ind_{C_a \wr S_b}^{S_{ab}} ) $ as follows.

\begin{Corollary}\label{cor:mashchar}
  Suppose $ \mash_a^b $ satisfies
  Properties (i) and (ii) in \Cref{prob:mash}. Then
    \[ \Ch( S^{\underline{\lambda}}\ind_{C_a \wr S_b}^{S_{ab}})
       = \sum_{\nu \vdash ab} \f{ \# \{ Q \in \SYT(\nu)
         : \mash_a^b(Q) = \ul \} }
         { \dim (S^{\underline{\lambda}} ) } s_{\nu}(\bfx), \]
  where $\mash_a^b(Q) \coloneqq \mash_a^b(w)$ for any
  $w \in \W_{ab}$ with $Q(w) = Q$.
\end{Corollary}

\begin{proof}
  We use, in order, \Cref{thm:grfrob_flexab}, Property
  (i), RSK, and Property (ii) to compute
  \begin{align*}
    \sum_{\ul \in \PP_a^b}
      \dim (S^{\ul})
      \Ch(S^{\ul} \ind_{C_a \wr S_b}^{S_{ab}}))
      q^{\underline{\lambda}}
      &= \W_{ab}^{\cont, \maj_a^b}(\mathbf{x}; q) \\
      &= \sum_{\alpha \vDash ab}
         \W_{\al}^{\maj_a^b}(q) \, \mathbf{x}^\alpha \\
      &= \sum_{\alpha \vDash ab}
         \W_{\al}^{\mash_a^b}(q) \, \mathbf{x}^\alpha \\
      &= \W_{ab}^{\cont,\mash_a^b}(\mathbf{x}; q) \\
      &= \sum_{\nu \vdash ab} (\SSYT(\nu)
         \times \SYT(\nu))^{\cont, \mash_a^b}(\mathbf{x}; q) \\
      &= \sum_{\nu \vdash ab} \SSYT(\nu)^{\cont}(\mathbf{x})
         \SYT(\nu)^{\mash_a^b}(q) \\
      &= \sum_{\nu \vdash ab} \SYT(\nu)^{\mash_a^b}(q)
         s_{\nu}(\mathbf{x}).
  \end{align*}
  The result follows by equating coefficients of
  $q^{\ul}$.
\end{proof}

\begin{Remark}
  When $a=1$ and $b=n$, we may replace $\ul$ with $\lambda \vdash n$.
  Under this identification, $\maj_1^n(w) = \sh(w)$,
  which clearly satisfies Properties (i) and (ii). When
  $a=n$ and $b=1$, we may replace $\ul$ with an element $r \in [n]$.
  Under this identification, we may set $\mash_n^1(w) = \maj_n(w)$,
  which satisfies Properties (i) and (ii). In this sense $\mash_a^b$
  interpolates between the major index $\maj_n$ and the shape under
  RSK, hence the name.
\end{Remark}

While $ \maj_a^b $ trivially satisfies Property (i),
it fails Property (ii) already when $a=b=2$, as in the following
example.
 
\begin{Example}\label{ex:mashnotmajab}
  Let $ v = 2314 $ and $ w = 1423 $. Then,
    \[ Q(v) = Q(w) = \Yvcentermath1{\young(124,3)} \]
  while
  \begin{align*}
    \maj_2^2(v) &= (\varnothing, (1, 1)) \\
    \maj_2^2(w) &= (\varnothing, (2)).
  \end{align*}
\end{Example}

\begin{Remark}
  When defining $\flex_a^b$ and $\maj_a^b$, we somewhat
  arbitrarily chose the lexicographic order on $\W_a$. Any other
  total order would work just as well. However,  $\maj_a^b$ continues to
  fail Property (ii) using any other total order when
  $a=b=2$ in \Cref{ex:mashnotmajab} since either $14<23$ or $23<14$.
\end{Remark}

\section*{Acknowledgements}

We would like to heartily thank our advisor, Sara Billey, for helpful
discussions and extensive comments on the manuscript. We would
also like to thank Vic Reiner, Christophe Reutenauer, Brendon Rhoades,
John Stembridge, and Shiela Sundaram for helpful conversations on this
material and related topics.

\bibliography{refs}{}
\bibliographystyle{alpha}

\end{document}